\documentclass[11pt]{amsart}
\usepackage{amssymb,amsmath,amsthm,amscd}
\usepackage{latexsym}
\usepackage[mathscr]{eucal}
\usepackage[english]{babel}
\usepackage{graphicx}
\newcommand{\so}{\mbox{${\mathfrak s \mathfrak o}$}}
\newcommand{\su}{\mbox{${\mathfrak s \mathfrak u}$}}

\newcommand{\un}{\mbox{${\mathfrak u}$}}

\newcommand{\af}{\mbox{${\mathfrak a}$}}

\newcommand{\ff}{\mbox{${\mathfrak f}$}}
\newcommand{\g}{\mbox{${\mathfrak g}$}}
\newcommand{\h}{\mbox{${\mathfrak h}$}}
\newcommand{\kf}{\mbox{${\mathfrak k}$}}

\newcommand{\m}{\mbox{${\mathfrak m}$}}
\newcommand{\n}{\mbox{${\mathfrak n}$}}
\newcommand{\p}{\mbox{${\mathfrak p}$}}

\newcommand{\s}{\mbox{${\mathfrak s}$}}

\newcommand{\C}{\mbox{${\mathbb C}$}}
\newcommand{\FF}{\mbox{${\mathbb F}$}}
\newcommand{\HH}{\mbox{${\mathbb H}$}}
\newcommand{\I}{\mbox{${\mathbb I}$}}

\newcommand{\PP}{\mbox{${\mathbb P}$}}

\newcommand{\R}{\mbox{${\mathbb R}$}}
\newcommand{\Z}{\mbox{${\mathbb Z}$}}

\newcommand{\tr}{{\rm tr}}
\newcommand{\ric}{{\rm Ric}}
\newcommand{\Ad}{{\rm Ad}}

\newcommand{\Aut}{{\rm Aut}}
\newcommand{\conv}{{\rm conv}}
\newcommand{\End}{{\rm End}}
\newcommand{\GL}{{\rm GL}}
\newcommand{\Ric}{{\rm Ric}}
\newcommand{\SO}{{\rm SO}}
\newcommand{\SU}{{\rm SU}}
\newcommand{\Sp}{{\rm Sp}}
\newcommand{\Spin}{{\rm Spin}}
\newcommand{\U}{{\rm U}}

\makeatletter
\def\numberwithin#1#2{\@ifundefined{c@#1}{\@nocnterrr}{%
  \@ifundefined{c@#2}{\@nocnterr}{%
  \@addtoreset{#1}{#2}%
  \toks@\expandafter\expandafter\expandafter{\csname the#1\endcsname}%
  \expandafter\xdef\csname the#1\endcsname
    {\expandafter\noexpand\csname the#2\endcsname
     .\the\toks@}}}}
\makeatother
\numberwithin{equation}{section}

\newtheorem{thm}[equation]{Theorem}

\newtheorem{prop}[equation]{Proposition}

\newtheorem{ex}[equation]{Example}
\newenvironment{example}{\begin{ex} \em}{\end{ex}}
\newtheorem{rem}[equation]{Remark}
\newenvironment{rmk}{\begin{rem} \em}{\end{rem}}

\begin{document}

\title{Einstein Metrics from Symmetry and Bundle Constructions: A Sequel}
\author{McKenzie Y. K. Wang}
\address{Department of Mathematics and Statistics,
McMaster University, Hamilton, Ontario, L8S 4K1, Canada}
\email{wang@mcmaster.ca}

\date{revised \today}

\begin{abstract} A survey was given by the author in 1999 \cite{Wa3} regarding developments
since Besse's volume \cite{Be} in the search for Einstein metrics via symmetry
reduction and bundle type constructions. The present article is a sequel to that
survey covering progress on selected topics since that time.
\end{abstract}

\maketitle


\bigskip
\setcounter{section}{-1}

\section{\bf Introduction}
It is a great honour for me to be asked to contribute an article to
commemorate the one hundredth anniversary of Professor Chern's birth.
This is especially so since I cannot claim any link via mathematical
genealogy to Professor Chern. While I was introduced to Chern classes
and Chern's proof of the Gauss-Bonnet theorem as an undergraduate,
I did not see him in person until I attended the Chern Symposium in 1979
as a graduate student. Later, in 1983, I met Professor Chern for the first time
at the MSRI, housed then in a building on Fulton Street in Berkeley.
I had the great fortune of being a member of the MSRI in the first year of
its operation. At that time I had just shifted my research from topology
to differential geometry, and the opportunity of participating in a special
year in differential geometry helped me enormously through this transition.
Even though I arrived after the Fall quarter, Professor Chern made me feel
at home immediately, and gave me sound advice and encouragement througout
my stay. I also met many geometers there for the first time, including Claude
LeBrun, with whom I edited a book years later, and Mario Micallef, through
whom I can claim to be a {\em collaborator} of a descendent of Professor Chern.

\medskip

This survey is concerned about the search for Einstein metrics and related
geometries through the use of symmetries and fibre bundle constructions.
We will focus on developments in the last twelve years, referring the
interested reader to \cite{Wa3}, which covers roughly the period 1987-1999,
and to \cite{Be} for even earlier work and foundations of the subject. The
present survey is not meant to be exhaustive, and will most certainly overlook
a number of important contributions. Our emphasis here will be on the case in
which the holonomy is generic, although some results about special holonomy metrics
will be discussed if they share a common method of construction with the generic
case. We also do not discuss Einstein warped products as these are best treated
in the context of quasi-Einstein metrics and the Bakry-Emery Ricci tensor. For these
topics we refer the readers to the articles \cite{Ca1}-\cite{Ca3}, \cite{HePW1},
\cite{HePW2}, \cite{WeWy} and the references therein.

We begin with the most symmetric situation, i.e., that of homogeneous Einstein
manifolds. Since it is known from \cite{AlK} that a Ricci-flat homogeneous manifold
is flat, our discussion is naturally divided into two cases corresponding to the
sign of the scalar curvature.  The condition of homogeneity reduces the Einstein
condition to a system of algebraic equations for which we seek {\em real} solutions
that satisfy an additional positivity condition. It is interesting to observe that for
both the positive and negative cases, the main conceptual advances in the last decade
result from applying a suitable variational approach. We shall deal mainly with the
positive case here (see \S2-\S3) as there is a superb survey of recent progress in the
negative case by J. Lauret \cite{La1}, who himself contributed immensely to the new
developments (see for example \cite{La2}-\cite{La4}).

We turn next to the case of metrics with {\em cohomogeneity one}. This means that
the isometry group of the metric acts with generic/principal orbits of codimension one.
The Einstein condition is now reduced to a system of non-linear ordinary differential
equations, together with additional conditions at finite and infinite boundary points
which ensure respectively smoothness and completeness of the metrics represented by the
solutions of the ODE system. While the general set-up and local issues for the cohomogeneity
one case are reasonably well-established (see \S 4), in contrast to the homogeneous case,
there is not yet a general theory for global existence. However, there have been some
efforts to find additional conserved quantities and interesting subsystems of the cohomogeneity
one Einstein equations (see \S 5), particularly in the Ricci-flat case. Cohomogeneity
one metrics with special holonomy arise in this way, and some progress on this front is
discussed in \S 5 and \S 6. One can also study the cohomogeneity one Einstein equations
using Painlev{\'e}-Kowalewski analysis (see \S 5). This is especially suited for exploring
integrability issues and for analysing the asymptotics of Einstein metrics.

Finally we survey recent efforts to construct complete Einstein metrics on fibre bundles
in \S 6. The approach here is based on a modification of the Kaluza-Klein ansatz. This
consists of finding suitable base and fibre metrics together with connections on the
associated principal bundle so that the resulting metric on the total space is Einstein.
In the case of a homogeneous base space and a cohomogeneity one fibre, this ansatz
reduces to the cohomogeneity one set-up. In almost all the non-cohomogeneity one cases
that have been studied to date, the fibre bundle has an abelian structural group. This
is because it is easy to satisfy the Yang Mills condition for the connection by using
Hodge theory. We believe that much remains to be understood for the case of non-abelian
structural groups.

While our survey deals mainly with Einstein metrics with generic holonomy, we have included
a discussion of complete cohomogeneity one metrics with holonomy ${\rm G}_2$ or $\Spin(7)$
in \S 6. We hope that this serves as an illustration of how the cohomogeneity one and bundle
viewpoints can interact in specific geometric situations.

I would like to thank Andrew Dancer and Spiro Karigiannis for their helpful suggestions and
remarks on preliminary versions of the paper. Partial support of my research through an Individual
Discovery Grant (OPG0009421) of the Natural Sciences and Engineering Research Council of
Canada is gratefully acknowledged.

\section{\bf Invariant Metrics on Compact Homogeneous Spaces}

We begin with some background information about invariant metrics on compact homogeneous
spaces which will be relevant throughout this survey.

\smallskip

A homogeneous manifold is a manifold which admits a transitive group of diffeomorphisms.
However, there can in general be infinitely many distinct transitive groups, i.e., non-conjugate
transitive subgroups of the diffeomorphism group, and these subgroups can be abstractly
isomorphic. This type of phenomenon was first discovered in \cite{WaZ2}. The simplest
example is given by the spaces $({\rm SU}(2) \times {\rm SU}(2))/{\rm U}(1)_{pq}$ where
$p, q$ are relatively prime integers. The underlying smooth manifold is $S^2 \times S^3$
by a classical result of Smale.

If we now fix a compact transitive Lie group $G$ on a homogeneous manifold $M$, then after choosing
a basepoint, we can write $M$ as the coset space $G/K$ where $K$ is the isotropy group
of the basepoint. In general, $G$ is neither the largest nor smallest Lie group (by inclusion)
that acts transitively on $M$. For example, the unit sphere $S^n$ is ${\rm SO}(n+1)/{\rm SO}(n)$,
but ${\rm SO}(n+1)$ lies in the full isometry group ${\rm O}(n+1)$, which is not connected,
as well as in the conformal group ${\rm CO}(n+1)$, which is not compact. When $n = 4m+3$,
then the subgroup ${\rm Sp}(m+1)$ still acts transitively on $S^n$, giving us the coset representation
${\rm Sp}(m+1)/{\rm Sp}(m)$. In going from a larger to a smaller transitive group, the space
of invariant metrics typically becomes larger. This increases the chances of finding Einstein
metrics, but the Einstein condition also becomes more complicated. From this viewpoint, the most
difficult case to analyse is that of group manifolds, although for a compact connected semisimple
Lie group the metric induced by the Killing form is well-known to be a non-negatively curved
positive Einstein metric. Indeed, as far as I know, the set of {\em all} left-invariant
Einstein metrics on any of the rank $2$ compact connected semisimple Lie groups have not been
completely determined.

We will consider connected homogeneous spaces $G/K$ where $G$ is a compact Lie
group, $K$ is a closed subgroup, and the $G$-action is {\em almost effective},
i.e., the kernel of the action is at most a finite group. Since we are interested
in positive Einstein metrics, in view of the theorem of Bonnet-Myers, we may as well
assume that the fundamental group is finite. {\em Unless otherwise stated, the above
will be our standing assumptions throughout \S 1 - \S 3.} We will also fix a
bi-invariant metric $b$ on $\g$, the Lie algebra of $G$, and use the normal metric
it induces on $G/K$, also denoted by $b$, as a background metric.

Using $b$ we obtain an ${\rm Ad}_K$-invariant orthogonal decomposition
$$ \g = \kf \oplus \p $$
where $\p$ is isomorphic to the tangent space of $G/K$ at the identity coset $[K]$. The ${\rm Ad}(K)$
action on $\p$ is called the {\em isotropy representation} of $G/K$. It is almost faithful when
the $G$-action is almost effective. The space of $G$-invariant Riemannian metrics on $G/K$ can be
identified with the space of all ${\rm Ad}_K$-invariant inner products on $\p$. It is a finite-dimensional
cone and hence is contractible. Its detailed description depends on the isotropy representation.

While we do not want to get into technical details here, we do want to give the reader
some feeling about the role played by the isotropy representation. In the case of a group
manifold, $K=\{1\}$, and so the isotropy representation consists of a direct sum of $\dim G$
trivial one-dimensional representations of $K$. The space of $G$-invariant metrics is just
the space of all inner products on $\p$, i.e., the symmetric space ${\rm GL}_{+}(d)/{\rm SO}(d)$
where $d=\dim G$. Note that there is always a basis of $\p$ with respect to which a given
$G$-invariant metric is diagonal, but of course there is no decomposition of $\p$ into
real irreducible summands such that the diagonal metrics with respect to that decomposition
include all the $G$-invariant metrics.

To parametrize the space of all $G$-invariant metrics, one can  decompose
$\p$ (using the background metric $b$) as an orthogonal direct sum
\begin{equation} \label{decomp}
 \p_0 \oplus \p_1 \oplus \cdots \oplus \p_{\ell}
\end{equation}
where the $\p_i$ are pairwise inequivalent ${\rm Ad}_K$ invariant summands, $\p_0$ is the
fixed point set of the isotropy representation, and each $\p_i$ is in turn a direct
sum of isomorphic irreducible {\em real} representations of ${\rm Ad}(K)$.
In general, the isotropy representation will have {\em multiplicities}, i.e., at least
one of the $\p_i$ consists of more than one irreducible summand. In this case, as
in the group manifold case, the further decomposition of $\p_i$ into irreducible
real subrepresentations is no longer unique (up to order). The number of parameters
of $G$-invariant metrics depends not only on the number of summands in $\p_i$ but
also on the nature of the summands.

\begin{example} \label{unitang}
Let $G={\rm SO}(m+2)$ and $K={\rm SO}(m)$ with $m \geq 3$. Then
$G/K$ is the unit tangent bundle of $S^{m+1}$. As $K \subset {\rm SO}(m+1) \subset G$
we see that the isotropy representation is $\I \oplus 2 \rho_m$, where $\I$ is the
trivial representation and $\rho_m$ is the vector representation of ${\rm SO}(m)$ on $\R^m$.
Notice that $\rho_m$ is {\em absolutely irreducible}, i.e., its complexification remains
irreducible over $\C$. The space of $G$-invariant metrics is $4$-dimensional, where the
additional parameter corresponds to the fact that the ${\rm SO}(m)$ equivariant
endomorphisms of $\R^m$ consist of multiples of the identity. When $m=2$, the representation
$\rho_2$ is no longer absolutely irreducible, and the space of $G$-invariant metrics
is $5$-dimensional (cf the next example).
\end{example}

\begin{example} \label{Aloff1}
Let $G={\rm SU}(3)$ and $K = {\rm U}_{1,-1}$ be the subgroup consisting of
diagonal matrices of the form ${\rm diag}(e^{i\theta}, e^{-i\theta}, 1)$.
The isotropy representation takes the form
$$  \I \oplus  V \oplus V \oplus W $$
where $\I$ is the trivial representation,  $V$ is the irreducible real
representation of ${\rm U}_{1, -1} \approx S^1$ corresponding to rotation by
$\theta$, and $W$ corresponds to rotation by $2\theta$. In the notation of
(\ref{decomp}), we have $\p_0 = \I, \, \p_1 = V \oplus V, \,\p_2 = W.$
The representation $V$ splits upon complexification into two irreducible unitary
representations which are dual to each other. Therefore, the $S^1$-equivariant
endomorphisms of $V$ form a $2$-dimensional space. As a result, the space of
$G$-invariant metrics is $1+2+2+1=6$-dimensional. Note that the subgroups
$\U_{1,0}$ and $\U_{0,1}$ are conjugate to $\U_{1,-1}$ in $\SU(3)$ and so
the corresponding homogeneous spaces are equivariantly diffeomorphic. Likewise
we can change the integers to their negatives and preserve equivariant
diffeomorphism.
\end{example}

\begin{example} \label{Aloff2}
Let $G={\rm SU}(3)$ and $K = {\rm U}_{1,1}$ be the subgroup consisting of
diagonal matrices of the form ${\rm diag}(e^{i\theta}, e^{i\theta}, e^{-2i\theta})$.
The isotropy representation takes the form
$$ 3\, \I \oplus V \oplus V$$
where $V$ is now the irreducible real representation of ${\rm U}_{1, 1} \approx S^1$
corresponding to rotation by $3\theta$, and $3\I$ corresponds to the subgroup
${\rm SU}(2)$ which commutes with ${\rm U}_{1,1}$. In the notation of (\ref{decomp}),
$\p_0 = 3\,\I$ and $\p_1 = V \oplus V$. As in the group manifold case, $\p_0$ contributes
$6$ parameters to the space of $G$-invariant metrics. As in the previous example,
the real irreducible summand $V$ splits upon complexification. We again pick up two
additional non-diagonal parameters of the $G$-invariant metrics to give
a total of $6+2+2=10$ parameters.
\end{example}

\begin{example} \label{quatsphere}
Let $G= {\rm Sp}(m+2)$ and $K={\rm Sp}(m)$ for $m \geq 1$. Using the inclusion
$K \subset {\rm Sp}(m)\times {\rm Sp}(2) \subset G$, we see that the isotropy representation
is
$$ 10\,\I \oplus [\nu_{m}]_{\R} \oplus [\nu_{m}]_{\R}$$
where $\nu_{m}$ is the $2m$-dimensional symplectic representation of ${\rm Sp}(m)$
and $[\nu_{m}]_{\R}$ denotes the corresponding irreducible real representation
on ${\R}^{4m}$. The latter representation splits upon complexification into
$2 \nu_{m}$. Because the space of ${\rm Sp}(m)$-invariant endomorphisms of
$[\nu_{m}]_{\R} \oplus [\nu_{m}]_{\R}$ is $4$-dimensional, the dimension of the
space of $G$-invariant metrics in this example is $55 + 2 + 4 = 61$.
\end{example}

The variational approach to the existence question for compact homogeneous spaces
is based on the Einstein-Hilbert action $\mathcal A$. This action associates to each metric
on a closed manifold the integral of its scalar curvature. The first variational formula
at a metric $g$ is given by
\begin{equation} \label{hilbert}
   d{\mathcal A}_g (h) = - \int_M \, g({\rm Ric}(g) - \frac{S_g}{2}g, h) \, d{\mu_g}
\end{equation}
where $S_g$ is the scalar curvature of $g$, $h$ is an arbitrary symmetric $2$-tensor,
and $d\mu_g$ is the Riemannian volume element.
If we take only constant volume variations, then the factor $2$ in the denominator
above should be replaced by the dimension $n$ of the manifold.

When $M=G/K$ and $g$ is a $G$-invariant metric, then ${\rm Ric}(g) - \frac{S_g}{n}g$ is
also a $G$-invariant symmetric $2$-tensor. So if $g$ is a critical point of the restriction
of $\mathcal A$ to the space of constant volume $G$-invariant metrics, then $g$ is also
a critical point of $\mathcal A$ on the space of all constant volume metrics, i.e., an
Einstein metric. (This argument is often referred to as the principle of symmetric
criticality.) It is often convenient to regard the metric g as a $b$-symmetric
automorphism $q$ of $\p$ via
\begin{equation} \label{gtoq}
 g(X, Y) = b(q(X), Y).
\end{equation}
We then have the {\em relative volume} $v(g)$ given by
\begin{equation} \label{relvol}
 d{\mu_g} = v(g)\,  d{\mu_b}, \,\,\, v(g)= \sqrt{{\det}\,g}.
\end{equation}
If we restrict ourselves to the space ${\mathscr M}_1^{G}(G/K)$ of volume $1,$ $G$-invariant
metrics on $G/K$, then the Einstein-Hilbert action gives rise to the
{\em scalar curvature function}
$$ S: {\mathscr M}_1^{G}(G/K) \rightarrow \R, \ \ S(g):={\mathcal A}(g)=S_g.$$
We will drop the reference to $G/K$ if no other space is being considered at the same
time.

Since $G$ is assumed to be compact, hence unimodular, $S_g$ is given by (see \cite{Be} (7.39))
\begin{equation} \label{scalar}
S_g = \frac{1}{2} \sum_{i,j} g^{ij} B(X_i, X_j)
     - \frac{1}{4} \sum_{i,j,k,l} \, g^{ik}g^{jl} g([X_i, X_j]_{\p}, [X_k, X_l]_{\p})
\end{equation}
where $B$ is the {\em negative} of the Killing form on $\g$, $\{X_i\}$ is any
$b$-orthonormal basis of $\p$, and $[\cdot, \cdot]_{\p}$ is the Lie bracket in $\g$
followed by $b$-orthogonal projection onto $\p$.
This formula shows that $S_g$ is a rational function on ${\mathscr M}_1^{G}$.
Notice that when $G$ is semisimple, we may let $b=B$; otherwise, $B$ is degenerate
on the centre of $\g$.

\begin{rmk} \label{recover}
Given a specific homogeneous space $G/K$, it is in practice simpler to compute
its scalar curvature function using the above formula than to compute the Ricci tensor
first and then take the trace. In situations where one needs to obtain the Ricci tensor,
it can be recovered from $S_g$ by differentiation as follows. Indeed, (\ref{hilbert}) leads to
$$ dS_g(h)\, v(g)\,{\rm vol}(b) + S_g \dot{v}\, {\rm vol}(b) = \frac{S_g}{2}\, g(g, h)\, v(g)\, {\rm vol}(b)
     -g({\rm Ric}(g), h) v(g)\, {\rm vol}(b),$$
where $h= \dot{g}$. But the second term on the left and the first term of the right
are equal since $v^{-1} \dot{v} = \frac{1}{2}\, \tr(q^{-1} \,\dot{q})=\frac{1}{2}\, g(g, h).$
(In the preceeding, $h(X, Y) = b(\dot{q}(X), Y).$) Hence
\begin{equation} \label{scal-to-ric}
      \tr(q^{-1}\hat{r}\,q^{-1}\dot{q}) = g({\rm Ric}(g), h)= -dS_g(h),
\end{equation}
where on the left-hand side $\hat{r}$ is the Ricci operator given by ${\rm Ric}(X, Y) = b(\hat{r}(X), Y).$
As this equation is linear in the components of $\hat{r}$ and holds for all variations $h$,
the Ricci tensor is completely determined.
\end{rmk}

\begin{rmk} \label{diagonal}
If we fix a $b$-orthogonal decomposition of $\p$ into real ${\rm Ad}_K$-irreducible
subrepresentations
\begin{equation} \label{decomp2}
 \p = \m_1 \oplus \cdots \oplus \m_r,
\end{equation}
then there is a useful formula for the scalar curvature of those $G$-invariant
metrics which are diagonal with respect to the above decomposition:
$$ g =  x_1 \, (b|_{\m_1}) \perp \cdots \perp x_r \, (b|_{\m_r}), \,\,\, x_i > 0, \,\,\,1 \leq i \leq r.$$
 This is
\begin{equation} \label{diagscal}
 S_g = \frac{1}{2} \sum_{i=1}^{r} \,\frac{{\beta}_i d_i}{x_i} - \frac{1}{4} \,\sum_{i,j.k} \,[ijk] \frac{x_k}{x_i x_j}
\end{equation}
where $\beta_i \geq 0$ are defined by $B|\m_i = \beta_i \, b|m_i$, $d_i := \dim_{\R} \m_i,$
and the non-negative coefficients
$$ [ijk] := \sum_{\alpha, \beta, \gamma} b( [X_{\alpha}, Y_{\beta}], Z_{\gamma})^2, $$
with $\{X_{\alpha}\}, \{Y_{\beta}\}, \{Z_{\gamma}\}$ being respectively $b$-orthonormal
bases of $\m_i, \m_j, \m_k$. Encoded in these coefficients is information about the
proper ${\rm Ad}_K$-invariant subalgebras of $\g$ which properly contain $\kf$.
We shall refer to these subalgebras as {\em intermediate subalgebras of the pair}
$(G, K)$.
\end{rmk}

Finally, we introduce the {\em gauge group} for the space ${\mathscr M}_1^{G}$.
This is the compact group $N(K)/K$ where $N(K)$ is the normalizer of $K$ in $G$.
An element $n \in N(K)$ acts on $G/K$ by sending $gK$ to $ngn^{-1}K$. Note that the
identity coset $[K]$ is a fixed point, and $N(K)$ acts on $\p$ via the adjoint representation.
The corresponding action on $G$-invariant metrics leads to isometric metrics. Hence
the scalar curvature function factors through ${\mathscr M}_1^{G}/(N(K)/K)$.
Unfortunately, this quotient space is still contractible, so one cannot apply Morse theory
in a simple-minded manner. In the following we illustrate the action of the gauge group
with some examples.

\begin{example} \label{Wallach1}
The gauge group for Example \ref{Aloff1} is the circle which is the quotient of
the subgroup of diagonal matrices in $\SU(3)$ by $\U_{1, -1}$.  A $G$-invariant metric
can be viewed as a positive definite symmetric block matrix of the form
$$ \left( \begin{array} {cccc}
        \rho & 0 & 0 & 0 \\
          0 & A & B & 0 \\
          0 & B^t & C & 0 \\
          0 & 0 & 0 & D
          \end{array} \right) $$
where $\rho$ is a positive real number, $A, C, D$ are positive multiples of the
$2 \times 2$ identity matrix, and $B$ is $2 \times 2$ depending on $2$ real parameters.
The gauge group action on $B$ can be viewed as rotation in the complex plane, and is
trivial on $A, C, D$. The space ${\mathscr M}_1^{G}/(N(K)/K)$ is therefore $5$-dimensional.

By \cite{Ni2}, there are two $G$-invariant Einstein metrics on
${\rm SU}(3)/{\rm U}_{1,-1}$. One of these is Sasakian Einstein, associated
with the circle fibration $\SU(3)/\U_{1, -1} \rightarrow \SU(3)/T$, where $T$ is the
usual maximal torus of $\SU(3)$. $B$ is $0$ for this metric, so it is fixed by the
gauge group. For the second Einstein metric, however, $B$ is nonzero. This fact will
become important in \S 6.
\end{example}

\begin{example} \label{Wallach2}
In Example \ref{Aloff2} above, the normalizer of ${\rm U}_{1,1}$ in ${\rm SU}(3)$
is the ${\rm SU}(2)$ embedded in the upper left-hand corner of ${\rm SU}(3)$. So
$N(K)/K \approx {\rm SO}(3)$. A $G$-invariant metric may be viewed as a positive
definite symmetric block matrix of the form
$$ \left( \begin{array} {ccc}
          A & 0 & 0 \\
          0 & B & C \\
          0 & C^t & D
          \end{array} \right) $$
where $A$ is $3 \times 3$, $B, D, C$ are $2 \times 2$, with $B, D$ being in addition
positive multiples of the identity. The action of the gauge group on $A$ is via conjugation.
Its action on the lower right $4 \times 4$ block can be thought of as the action
of ${\rm SO}(3)$ on the quaternions, where the real axis consists of multiples of the
$4 \times 4$ identity matrix and ${\rm SO}(3)$ acts in the usual way on the imaginary
quaternions. So a $G$-invariant metric is isometric either to one where $A$ is diagonal
or one where $C=0$. The space ${\mathscr M}_1^{G}/(N(K)/K)$ is $7$-dimensional.

By \cite{Ni2} there are, up to isometry, two $G$-invariant Einstein metrics
on ${\rm SU}(3)/{\rm U}_{1,1}$. Both metrics are associated with the Riemannian submersion
${\rm SO}(3) \rightarrow  {\rm SU}(3)/{\rm U}_{1,1} \rightarrow \C\PP^2$, where the base
is viewed as a self-dual manifold (cf Proposition 14.85 in \cite{Be}). Hence, in terms of
the above block form, we have $C=0$, and $A$ and $B=D$ are  multiples of the identity.
So the Einstein metrics are fixed points of the gauge group.
\end{example}

\begin{example} \label{spin4}
Let us take $G={\rm SU}(2) \times {\rm SU}(2)$ and $K= \{(1, 1)\}$, the trivial subgroup.
A left-invariant metric on $\g$ may be viewed as a positive definite symmetric block matrix
$$ \left( \begin{array}{cc}
        A  &  B  \\
        B^t & C
   \end{array}  \right)$$
where $A, B, C$ are $3 \times 3$ sub-matrices.  Now $(a_1, a_2) \in N(K) = G$ acts on
a metric as above by conjugation by
$$\left( \begin{array}{cc}
              P_1  &  0 \\
              0     &  P_2
   \end{array}  \right) $$
where $P_i$ are respectively the images of $a_i$ in ${\rm SO}(3)$ under the usual covering map.
Thus one sees that an arbitrary left-invariant metric is isometric either to one where
$A$ and $C$ are diagonal or to one where $B$ is diagonal. The space ${\mathscr M}_1^{G}/(N(K)/K)$
is $15$-dimensional, so the Einstein condition consists of $15$ equations in $15$ unknowns.

Two left-invariant Einstein structures have been found on $G$. One is given by the product metric,
which corresponds to $A=C=\lambda I$ with $\lambda > 0$ and $B=0$. It is a fixed point under
the gauge group. The second known Einstein structure (cf \cite{Jen1}) is given by the Killing form metric
on $(G \times G \times G)/\Delta G$ where $\Delta G$ is the diagonally embedded subgroup.
It corresponds to $A=B=\frac{2}{3}I, B=-\frac{1}{3}I$. The orbit of this metric under the
action of the gauge group is ${\rm SO}(3)$. It has been shown in \cite{NiR} that there are no
further left-invariant Einstein metrics which have an addition circle symmetry.
\end{example}

\begin{rmk}
The fact that in Example \ref{spin4} there is an Einstein metric whose orbit under the gauge
group is $\SO(3)$ shows that the Einstein equations are in general non-generic in the sense
of algebraic geometry. Therefore, even though the assumption of homogeneity reduces the
Einstein condition to a system of algebraic equations, techniques based on analytic concepts
may really be indispensable in their study.
\end{rmk}

\section{\bf Variational Approach for Positive Homogeneous Einstein Metrics}

The variational approach to finding homogeneous Einstein metrics was first introduced by
G. Jensen \cite{Jen1}. It was then systematically developed in \cite{WaZ1}, \cite{BWZ},
and  \cite{Bo1}. In its essence, the variational approach consists of two interdependent
parts:
(i) investigating the analytic properties of the scalar curvature function $S_g$, and
(ii) finding (easily) computable combinatorial invariants of the pair $(G, K)$ which
would guarantee the existence of critical points of $S$.

The first analytic property that comes to mind is boundedness of $S$. The formula
(\ref{diagscal}) suggests that $S$ should be negative in most regions. However, this
first intuition has to be modified by the following observations.

\begin{itemize}
\item Let $\h$ be an intermediate subalgebra whose corresponding Lie group $H$ is compact.
    If we shrink the metric $b$ along the fibres of the fibration
    $H/K \rightarrow G/K \rightarrow G/H$ (while keeping the volume constant), then
    the scalar curvatures of the resulting metrics tend to $+\infty$ if $H/K$ is 
    non-abelian and they tend to $0$ if $H/K$ is a torus.
\item If $g_t$ is a curve of $G$-invariant metrics starting from the basepoint $b$ and
   going ``radially" to infinity, then the scalar curvatures along it are bounded from below
   only when ``many" of the coefficients $[ijk]$ vanish.
\item Therefore, for ``most" $G/K$, we cannot expect the scalar curvature function to be
  bounded from above or from below. So the critical points of $S$ on ${\mathscr M}_1^{G}$
  are generally saddle points.
\end{itemize}

Reflecting further on the above intuitive picture raises some technical issues.

\begin{enumerate}
\item[a.] For a connected $G/K$ with $G$ compact, the identity component $G_0$ of $G$
   still acts transitively with isotropy group $G_0 \cap K$. We may as well assume $G$ to
   be connected since the set of $G_0$-invariant metrics contains the set of
   $G$-invariant metrics. Assuming $K$ to be connected is a different matter because
   this corresponds to looking at invariant Einstein metrics on the cover $G/K_0$, which
   is a different manifold. One then has the problem of deciding which of the Einstein
   metrics found there descend to $G/K$. For a qualitative study of the Einstein condition,
   separating the problem in this way is undesirable. This explains the assumptions made
   in \cite{BWZ} and \cite{Bo1}, and some of the technicalities in these papers.
\item[b.] Since $K$ is not connected in general, not every subalgebra $\h$ satisfying
    $\kf \subset \h \subset \g$ is $\Ad_{K}$-invariant. This is why $\Ad_K$-invariance
    is included in the notion of an intermediate subalgebra of $(G, K)$. Furthermore,
    the analytic subgroup $H_0$ with $\h$ as Lie algebra contains $K_0$ but not necessarily
    $K$. So one has to consider $H$ which is the subgroup generated by  $H_0$ and $K$.
    However, $H$ need not be closed as $H_0$ need not be closed. This occurs when
    $\h = \kf \oplus \af$ where $\af$ is an abelian subalgebra in $\p$ that is irrationally
    embedded in a strictly larger compact abelian subalgebra. Indeed, such an $\h$ lies
    in a continuous family of subalgebras of the same dimension.
\item[c.] The formula (\ref{diagscal}) only applies for invariant metrics which are diagonal
   with respect to a fixed $b$-orthogonal decomposition (cf (\ref{decomp2})) of $\p$
    into $\Ad_K$-invariant, real irreducible summands. However, every invariant metric
    is of diagonal form with respect to at least one decomposition of $\p$, and the space
    of all decompositions as in (\ref{decomp}) of $\p$ is compact. In obtaining estimates
    of the scalar curvature function $S$, one has to work with Eq. (\ref{diagscal}) in such
    a way that all constants appearing in the estimates are independent of the choice
    of decomposition.
\end{enumerate}

\smallskip
\noindent{\bf 2A. Compactness Properties}
\smallskip

We will now describe some specific results in \cite{WaZ1}, \cite{BWZ} and \cite{Bo1}.
We remind the reader of our standing assumptions (see \S 1), which may be suppressed from
the statements of theorems below. The first existence result based on bounding the
scalar curvature is

\begin{thm} $($\cite{WaZ1}, \cite{BWZ}$)$ \label{primitive}
Let $G/K$ be a primitive homogeneous space, i.e., one whose only $G$-invariant foliations
are the ones by the whole manifold or by its points. Then its scalar curvature function is
bounded from above and proper. In particular, it has a global maximum, which is a $G$-invariant 
Einstein metric.
\end{thm}

\begin{rmk} Although this result has been generalized in \cite{Bo1}, we state it here
for various reasons. First, the proof of the above theorem is essentially the
same as that of Theorem 2.2 in \cite{WaZ1}. In particular, it does {\em not} require
the Palais-Smale condition, to be described shortly. It is the simplest example of how
the scalar curvature function may be estimated, and is a good starting point for readers
who may be interested in the more complicated proofs of the theorems below. Second,
when $K$ is connected the converse also holds, and so in this case one obtains an
analytic characterization of primitivity.

Recall that the classification of primitive homogeneous spaces (without compactness
assumptions) is a problem posed originally by S. Lie. It includes
the classification of maximal subalgebras in real Lie algebras, but in its global form
also includes the situation when $G$ and $K$ are not necessarily connected. Partial
classifications have been given by Lie himself, Morozov, Dynkin, and M. Golubitsky and
B. Rothschild \cite{Go}, \cite{GoR}. Komrakov \cite{Kom} gave a classification
assuming $G$ is connected. In unpublished work, Ziller and I have proved that a
connected effective homogeneous space $G/K$ with $G$ compact (but not necessarily
connected), $K$ closed, and ${\rm rank}(G) = {\rm rank}(K)$ is primitive iff $G/K$
has irreducible isotropy representation.
\end{rmk}

In order to use the variational method to find saddle points of $S$, we need the following
condition to hold.

\smallskip

\noindent{\bf Palais-Smale Condition:} Suppose $g_i$ is a sequence in ${\mathscr M}_1^{G}$
such that $S(g_i)$ is bounded and  $\|{\rm grad}S \|$ converges to $0$, where $\| \cdot \|$
is the $L^2$ norm restricted to ${\mathscr M}_1^{G}$. Then there is a  subsequence that converges
to an element of  ${\mathscr M}_1^{G}$ in the $C^{\infty}$ topology.

\smallskip

\begin{thm} $($\cite{BWZ}, Theorem 1.1$)$ \label{PS} For every $a> 0$ the above Palais-Smale
condition is satisfied by $S$ on the subset $\{g \in {\mathscr M}_1^{G}: S_g \geq a\}$.
\end{thm}

\begin{rmk}
There is currently no elementary proof of the above theorem except in the case where
rank $G =$ rank $K$. Even in this case the proof is rather unwieldy. Instead,
one appeals to the deep work of Cheeger-Colding (cf \cite{ChCo}) on Gromov-Hausdorff
convergence of Riemannian manifolds under lower Ricci bounds. The point is that the assumptions
in the Palais-Smale condition give a lower Ricci bound, which in turn implies, by the Bishop-Gromov
volume comparison theorem, a uniform lower bound on the volumes of unit metric balls. Hence the
Cheeger-Colding theory can be applied, and the sequence of homogeneous Riemannian manifolds
subconverges in the pointed Gromov-Hausdorff topology to a pointed complete metric space of
the same dimension. Normally, one only has $C^{1, \alpha}$ convergence of the metrics to a
metric on an open subset of the limit space. However, using homogeneity one in fact obtains
$C^{1, \alpha}$ subconvergence to a smooth Einstein manifold, {\em modulo gauge transformations}.

In order to have the Palais-Smale condition as stated above, we must still keep track of
the original transitive $G$-actions and the subconvergence of the original sequence of metrics.
This can be done by using an equivariant version of the Gromov compactness theorem due to
Fukaya \cite{Fu} and some classical results of Montgomery-Zippin \cite{MoZ}. It turns out that
the limit Einstein manifold has a limit transitive $G$-action such that the isotropy group of the
basepoint is conjugate to $K$. This in turn implies the subconvergence of $g_i$ in the
smooth topology to a $G$-invariant Einstein metric.
\end{rmk}

\smallskip
\noindent{\bf 2B. Moduli}
\smallskip

The Palais-Smale condition immediately allows one to deduce some properties of the
moduli space of homogeneous Einstein metrics.
Let ${\mathscr E}(G/K) \subset {\mathscr M}_1^{G}$ be the subset of Einstein metrics
lying in the set of volume $1$ $G$-invariant metrics on $G/K$.
It is a semialgebraic set, since it is given by the Einstein equations
and some inequalities expressing the positive-definiteness of the solutions
as symmetric $2$-tensors.  Hence it has a local stratification by real analytic
submanifolds, is locally path connected, and has finitely many topological components.
Together with Theorem \ref{PS}, one obtains

\begin{thm} $($\cite{BWZ}, Theorem 1.6$)$ \label{compactmoduli}
${\mathscr E}(G/K)$ has finitely many components, each of which is compact.
$S$ takes on only finitely many values on ${\mathscr E}(G/K)$.
\end{thm}

Let us turn next to the situation where $M$ is a fixed compact connected
homogeneous manifold {\em and} we allow the transitive actions to change.
Let  $(M=G_i/K_i, g_i)$ be a sequence of unit volume $G_i$-invariant Riemannian
metrics with scalar curvatures $S_{g_i}$ converging to $a>0$ and
$\|\ric^0(g_i)\|_{g_i}$ converging to $0$. By replacing $G_i$ by its identity component
and noting that only finitely many abstract isomorphism classes of compact connected Lie
groups can act transitively and almost effectively on $M$, one can assume that in the
above sequence all $G_i$ are the same compact connected Lie group $G$. Of course,
the transitive actions will in general be different, i.e., the isotropy groups
need not be conjugate in $G$.  But the proof of the Palais-Smale condition
now shows that a subsequence of $g_i$ converges in the $C^{\infty}$ topology to a
$G$-homogeneous Einstein metric on $M$. Using this one deduces

\begin{thm} $($\cite{BWZ}, Theorem 1.8$)$ \label{homogmoduli}
Let $M$ be a compact, connected homogeneous manifold and ${\mathscr E}_h(M)$ denote
the subspace of homogeneous Einstein metrics on $M$ in ${\mathscr M}_1$, the space
of all unit volume smooth Riemannian metrics on $M$ equipped with the $C^{\infty}$ topology.
For any positive number $a$, the subspace of ${\mathscr E}_h(M)/{\rm Diff}(M)$ consisting
of all Riemannian structures with scalar curvature $\geq a$ has finitely many components
and each component is compact. Furthermore, only finitely many values of the scalar curvature
are assumed in this space.
\end{thm}

\begin{rmk} The above theorem implies that on a compact connected homogeneous
Einstein manifold there is a maximum value for the scalar curvature among all
homogeneous Einstein metrics of unit volume.

On the other hand, the example of $M=S^2 \times S^3 = (\SU(2) \times \SU(2))/\U_{p,q}$
(with $p, q$ relatively prime) mentioned at the beginning of \S 1 actually
has a unique unit volume $\SU(2) \times \SU(2)$-invariant Einstein metric whenever
$(p,q) \neq (1,1), (1, 0)$ or $(0, 1)$ (cf \cite{WaZ2}). The sequence of Einstein
constants has $0$ as the only limit point.
\end{rmk}

\vspace{2cm}

\noindent{\bf 2C. Graph Theorem}
\smallskip

The Palais-Smale condition also allows one to apply the gradient flow of $S$ on
${\mathscr M}_1^G$ (which coincides with the restriction of the normalized Ricci flow)
to obtain critical points by studying the structure of {\em superlevel sets} of $S$, i.e.,
sets of the form $\{g \in {\mathscr M}_1^G: S(g) \geq a \}.$

For an intuitive account of how such a study can be undertaken, one can imagine
${\mathscr M}_1^G$ to be an open ball with the background metric $b$ as its centre.
Since $b$ has positive scalar curvature, $S$ is positive in some neighbourhood of $b$.
Let $H$ be a closed subgroup such that $K \subset H \subset G$ and $H/K$ is not a
torus. Then we can shrink $b$ along $H/K$ and make $S$ increase to $+\infty$ along
a ray from $b$ to the boundary of the ball. Next choose a sufficiently large sphere
$\Sigma$ enclosing $b$. Rays of the type just described could form continuous
families, but their intersections with the sphere $\Sigma$ should lie in a finite number
of pairwise disjoint open subsets of the sphere whose union is {\em not} all of
$\Sigma$. Assume there are at least two such open sets, and let $\Omega$ be the region
of the ball that remains after we remove the closed region bounded by $\Sigma$
together with all the rays passing through the finitely many open subsets we singled
out in $\Sigma$. If we can show furthermore that $S$ is negative on $\Omega$,
then there must be a critical point of $S$ in the closed ball bounded by $\Sigma$.
Otherwise, we can take the union of two rays emanating from $b$ and passing through
two of the disjoint open sets in $\Sigma$. This curve lies in some superlevel set of
$S$ corresponding to a positive level. By the Palais-Smale property, this curve, which
is anchored at infinity (as $S$ tends to $+\infty$ there), must leave the region
bounded by $\Sigma$ in finite time. We therefore obtain a curve which must pass through
a negative scalar curvature region and still connects two different components of a
high superlevel set of $S$. This contradicts the fact that the gradient flow increases $S$.

To make the above statements rigorous, we introduce the {\em graph} of $G/K$.
This is a graph whose vertices are essentially the $\Ad_K$-invariant intermediate
subalgebras of the pair $(G, K)$ mentioned at the end of Remark \ref{diagonal}.
However, since these intermediate subalgebras can come in continuous families,
we actually define a {\em vertex} $[\h]$ to be the connected component of $\h$ in the set of
all intermediate subalgebras of the same dimension, regarded as a subset of
the Grassmannian of $\dim \h$-subspaces in $\g$. Two vertices are connected by an edge
when a subalgebra in one vertex is contained in or contains a subalgebra in the other
vertex. Properties of semialgebraic sets allow one to deduce that the graph has only
finitely many vertices.

A connected component of the above graph is called a {\em toral component} if {\em all}
intermediate subalgebras $\h$ occurring in the component have the property that
$\h/\kf$ are abelian. A component that is not toral is called {\em non-toral}.
In Example \ref{unitang}, the graph of $(\SO(m+2), \SO(m)), m \geq 3$
consists of two disjoint vertices $[\so(m+1)]$ and $[\so(m)\oplus \so(2)]$. The second
vertex gives a toral component. When $m=2$, the graph consists of the vertex
$[\so(3)]$ and a second non-toral component containing the vertices $[\so(2)\oplus\so(2)],
[\su(2)], [\su(2)^*], [\un(2)],$ and $[\un(2)^*]$. (Note that there are two conjugacy classes of
$\un(m)$ in $\so(2m)$ for $m \geq 2$.)

\begin{thm} $($Graph Theorem, \cite{BWZ}, Theorem 3.3$)$ \label{graphthm}
Suppose the graph of $G/K$ has at least two non-toral components. Then $G/K$ admits
a $G$-invariant Einstein metric. The same holds if $G$ and $K$ are both connected
and the graph of $G/K$ has at least two components.
\end{thm}

\begin{rmk} \label{graphrmk}
(i) The graph theorem is ineffective when $G$ is non-semisimple, as in this case
the graph is actually always non-empty and connected (cf Proposition 4.9 in \cite{BWZ}).
However, under our standing assumption of finite fundamental group, the semisimple part
of $G$ always acts transitively on $G/K$.

(ii) The critical points of $S$ produced by the graph theorem have coindex at most $1$.
Recall that the coindex is the number of positive eigenvalues of the Hessian at a critical
point.
\end{rmk}

The graph theorem shows immediately that $\SO(m+2)/\SO(m)$ has an invariant Einstein
metric. (This metric was first discovered by S. Kobayashi \cite{Ko}, and is actually
the only invariant one  when $m\geq 3$ \cite{BaHs}.) Further examples of the use of
the graph theorem can be found in \cite{BWZ}. We mention here one example from
\cite{BWZ} which shows the importance of not assuming $G$ and $K$ to be connected in
developing the general theory.

\begin{example}
Let $G_0 = \SO(k(p+q)), K_0 = \SO(p)^k \times \SO(q)^k$ with $p, q \geq 3, k\geq 2.$
Then the graph of $G_0/K_0$ is actually connected. However, if we let $\Gamma$ to be
the symmetric group on $k$ letters, then $\Gamma$ can be made to act on $G_0$ such that
it also permutes the $\SO(p)$ and the $\SO(q)$ factors among themselves. Let
$G = G_0 \ltimes \Gamma, H= H_0 \ltimes \Gamma$ be respectively the semidirect products.
Then the graph of $G/K$ has at least two non-toral components. These are the components
which contain the vertices $[\so(kp) \oplus \so(kq)]$ and $[k\so(p+q)]$. Therefore, the
existence of an invariant Einstein metric comes about by considering disconnected transitive
groups.
\end{example}

\medskip

\noindent{\bf 2D. Higher Coindex Critical Points}

\smallskip

In order to produce critical points of higher coindex, one has to look even more closely
at the analytic properties of $S$ and how these interact with the topology of the high
superlevel sets of $S$. This has been achieved by C. B\"ohm in \cite{Bo1}.
In this work, a simplicial complex is constructed from the intermediate subalgebras of
the pair $(G,K)$, as in the case of the graph of $G/K$. Non-contractibility of the
simplicial complex implies the existence of critical points of $S$ of higher coindex.

To describe the simplicial complex, we consider the restriction of the isotropy
representation to the identity component $K_0$ of $K$. Denote by
$\ff \subset \p$ the subspace on which $K_0$ acts trivially. As is easily seen,
$[\ff, \ff] \subset \ff$, and $\kf \oplus \ff$ is the Lie algebra of the identity
component of $N_G(K_0)$, the normalizer of $K_0$ in $G$. One now fixes a maximal
torus $T$ in the compact connected Lie group $N_G(K_0)^0/K_0$, whose Lie algebra
is $\ff$. Consider the set of intermediate subalgebras of $(G, K)$
which are in addition $\Ad_T$-invariant. An $\Ad_T$-invariant intermediate subalgebra
$\h$ is {\em minimal non-toral} if it is non-toral (i.e., $\h/\kf$ is non-abelian) and
any $\Ad_T$-invariant intermediate subalgebra $\h^{\prime}$ contained properly in $\h$
must be toral (i.e., $\h^{\prime}/\kf$ is abelian).

It turns out that there can only be finitely many minimal non-toral $\Ad_T$-invariant
intermediate subalgebras. These subalgebras generate a finite number of $\Ad_T$ invariant
(non-toral) intermediate subalgebras of $(G, K)$. Flags of these intermediate subalgebras
give rise to a finite simplicial complex in the usual way: the intermediate subalgebras
themselves are the vertices, a pair of intermediate subalgebras $\h_1 \subset \h_2$
forms a $1$-simplex, a triple of intermediate subalgebras $\h_1 \subset \h_2 \subset \h_3$
forms a $2$-simplex etc.

The main existence theorem is then

\begin{thm} $($\cite{Bo1}, Theorem 1.5$)$ \label{simcplxthm}
Let $G/K$ be a compact connected homogeneous space with $G$ compact Lie, $K$ closed.
Let $T$ be a maximal torus in $N_G(K_0)^0/K_0$ and $\Delta^T(G/K)$ be the simplicial
complex defined as above. If it is not contractible, then $G/K$ admits a $G$-invariant
Einstein metric.
\end{thm}

\begin{rmk} \label{simp-props}
(i) In general, the simplicial complex $\Delta^T(G/K)$ depends on the choice of $T$. However,
this is not the case when $G$ and $K$ are both connected. In this situation, there is
a one-to-one correspondence between the minimal non-toral $\Ad_T$-invariant intermediate
subalgebras of $(G, K)$ and the minimal non-toral intermediate subalgebras of
$(G, T\cdot K)$, where $T\cdot K$ denotes the subgroup of $G$ generated by $T$ and $K$,
which is compact and connected. The gauge group of $G/(T\cdot K)$ is zero-dimensional.
By Proposition 4.1 and Corollary 4.5 in \cite{BWZ} there can only be finitely many
intermediate subalgebras for $(G, T \cdot K)$. This explains why the simplicial complex
has finitely many vertices. The general case follows from this using the finiteness of
the fundamental group.

(ii) The Einstein metric in the above theorem corresponds to a critical point of
$ S \left|{\mathcal M}_1^G \right.$. In \cite{Bo1} it is further shown that if there is a field $\FF$
and a non-negative integer $q$ such that $\tilde{H}_q(\Delta^T(G/K); \FF) \neq 0$, then
the maximal subspace on which the Hessian of $S \left|{\mathcal M}_1^G \right.$ (at the critical point)
is positive semidefinite has dimension $\geq q+1$.
\end{rmk}

For concrete applications as well as theoretical reasons, B\"ohm introduced some variants
of the simplicial complex defined above:

\smallskip

I. A homogeneous space of the type under consideration is said to be of {\em finite type}
if there are only finitely many minimal non-toral intermediate subalgebras for the pair
$(G, K)$. In this case one can consider the set of intermediate subalgebras (necessarily
non-toral and finite in number) generated by these minimal ones. One then obtains
a simplicial complex $\Delta^{\rm min}(G/K)$ whose $k$-simplices correspond to flags made up
of $k$ intermediate subalgebras $\h_1 \subset \cdots \subset \h_k$ from the afore-mentioned
finite set.

Examples of homogeneous spaces of finite type include $G/K$ whose isotropy representation
has no multiplicities (i.e., the $\Ad_K$ irreducible subrepresentations are pairwise
non-isomorphic), and $G/K$ for which $N_G(K_0)^0/K_0 = \{1\}.$ On the other hand, in
Example \ref{unitang}, $N_G(K_0)^0/K_0 \approx \SO(2),$ thereby giving rise to a circle of
non-toral intermediate subalgebras of type $\so(m+1)$. Hence $\SO(m+2)/\SO(m)$ is not of
finite type. Likewise, Example \ref{spin4} is not of finite type since there is an $\R\PP^2$
family of non-toral intermediate subalgebras of type $\su(2)\oplus \un(1)$.

Non-contractibility of $\Delta^{\rm min}(G/K)$ also implies the existence of a $G$-invariant
Einstein metric (cf Theorem 1.4, \cite{Bo1}). In fact, the main case of the proof for
Theorem \ref{simcplxthm} is that when $\n(\kf)=\kf$, which is a special case of the finite
type condition.

\smallskip

II. One often considers $G/K$ where $G$ and $K$ are both connected. (Indeed, $G/K$ is finitely
covered by $G_0/K_0$, where $G_0$ and $K_0$ are the respective components of the identity.)
In this situation, the simplicial complex $\Delta^T(G/K)$ is independent of the choice of
the maximal torus $T \subset N_G(K)^0/K$. Furthermore, since $\pi_1(G/K)$ is assumed to be
finite, Remark \ref{simp-props}(i) implies that $\Delta^T(G/K) \approx \Delta(G/TK) =
\Delta^{\rm min}(G/TK)$. Hence one only needs to examine the non-contractibility of this
last simplicial complex.

For example, let $G=\SU(3)$ and $K=\U_{p,q}$ with $p, q$ relatively prime for convenience.
$G/K$ are the simply connected Aloff-Wallach spaces. The gauge group $N(K)/K$ is $1$-dimensional
except when $(p, q)=\pm(1,1), \pm(1, -2), \pm(2, -1)$, in which case it is $\SO(3)$,
cf Example \ref{Wallach2}. Therefore, we can take $G/(T \cdot K)$ to be $\SU(3)/T^2$ in all cases.
The simplicial complex of $\SU(3)/T^2$ consists of $3$ vertices, which correspond to the
three subalgebras of type $\su(2)$ determined by the three pairs of roots of $\su(3)$.
Hence one gets at least one $SU(3)$-invariant Einstein metric on each Aloff-Wallach space
without any computations.

\smallskip

The example just given is in fact a special case of an effective way of producing
compact homogeneous spaces with non-contractible simplicial complexes described in
\cite{Bo1}. Suppose for simplicity that $G/K$ is in addition simply connected with $G$
connected. As the semisimple part of $G$ still acts transitively on $G/K$, we may
assume further that $G$ is simply connected and semisimple. These assumptions easily
imply that $G/K$ is a homogeneous torus bundle, possibly with zero dimensional
fibres, over a product of spaces $G_i/K_i$ of the same type (i.e., with $G_i$ compact,
simply connected, semisimple, $K_i$ closed connected) having the additional property
that the gauge groups are zero-dimensional. In \cite{Bo1} such homogeneous spaces
are called the {\em prime factors} of $G/K$.

Using the arguments in II above and certain well-known topological facts due to
Quillen \cite{Q} and Milnor (homology of joins) \cite{Mi}, one obtains

\begin{thm} $($\cite{Bo1}, Theorem B$)$ \label{noncontract}
Let $G/K$ be a homogeneous space where $G$ is compact, simply connected, semisimple
and $K$ is a closed connected subgroup. If the simplicial complexes of its prime
factors all have non-trivial reduced homology $($wrt some field of coefficients$)$,
then the simplicial complex of $G/K$ is not contractible, and hence $G/K$
admits a $G$-invariant Einstein metric.
\end{thm}

Examples of prime homogeneous spaces with non-contractible simplicial complexes
given in \cite{Bo1} include:

\begin{itemize}
\item $G/T$ where $G$ is a simple classical Lie group and $T$ a maximal torus,
\item ${\rm G}_2/\U(2), {\rm F}_4/(\U(3)\cdot \Sp(1)), {\rm E}_6/(\Spin(10)\cdot \SO(2)),
    {\rm E_7}/({\rm E}_6 \cdot \SO(2)), {\rm E}_8/({\rm E}_6 \cdot \U(2))$,
\item $(\SU(n_1) \times \cdots \SU(n_k))/K$ where $\kf$ is a regular subalgebra of $\g$
      in the sense of Dynkin. This means that $\kf$ is obtained from $\g$ by choosing
      an abelian subalgebra of a maximal abelian subalgebra of $\g$ together with a
      subset of the root spaces of $\g$,
\item$\SO(n_1 + \cdots + n_k)/(\SO(n_1) \times \cdots \times \SO(n_k))$ where $k\geq 2$
      and $\n_i \geq 2$.
\end{itemize}
Using these prime homogeneous spaces and Theroem \ref{noncontract}, one can list
numerous compact simply connected homogeneous spaces admitting invariant Einstein
metrics.

\smallskip

The more refined scalar curvature estimates necessary for the proof of Theorem
\ref{simcplxthm} also provide a generalization of Theorem \ref{primitive}.

\begin{thm} $($\cite{Bo1}, Theorem 5.22$)$ \label{upperbound}
Let $G/K$ be a compact connected homogeneous space with $G$ compact Lie, $K$ closed.
The scalar curvature function is bounded from above iff $(G, K)$ has no non-toral
$\Ad_K$-invariant intermediate subalgebras. In this case it has a global maximum
which corresponds to a $G$-invariant Einstein metric.
\end{thm}

\section{\bf  Positive Homogeneous Einstein Manifolds: Classifications and Non-existence}

In contrast to general existence theorems, classification results require an
enumeration of the relevant homogeneous spaces $G/K$, as well as finding {\em all}
solutions to the corresponding Einstein equations. Finally, one has to distinguish
the resulting Einstein metrics up to isometry and homothety, which could prove
difficult. If one is not exhausted already, there is still the problem of classifying
the coset spaces $G/K$ up to diffeomorphism.

Regarding the enumeration problem, given a homogeneous manifold $M$, one first has to
determine all compact Lie groups $G$ which act transitively on $M$. There can in general be
infinitely many such groups which are not conjugate to each other in ${\rm Diff}(M)$.
Next suppose we have a fixed coset $G/K \approx M$. Since $M$ is assumed to be connected,
the identity component of $G$ still acts transitively on $M$. So we may assume $G$ to be connected.
This leaves open the question of determining the full isometry group of any $G$-homogeneous
metric we may find. If there is a subgroup of $G$ that is still transitive on $M$, then
invariant metrics on the coset space of the smaller transitive group include those
of $G/K$. Since the fundamental group of $G/K$ is assumed to be finite, the semisimple
part of $G$ is a transitive subgroup. So we may assume $G$ to be semisimple. Since the
group action is assumed to be almost effective rather than effective, we may further
suppose that $G$ is simply connected. In this situation, if $K_0$ is the identity component
of $K$, then $G/K_0$ is the universal cover of $M$. Once one finds all the $G$-invariant
Einstein metrics on $G/K_0$, the ones which factor through the action of $K/K_0$ are the
$G$-invariant Einstein metrics on $G/K$. These remarks explain the hypotheses in most
of the classification theorems we describe below.

\smallskip

{\bf (a) Six-dimensions:} (\cite{NiR}) The classification here assumes that $G$ is
     compact, connected, semisimple and $G/K$ is simply connected. Besides
     symmetric spaces, the other possibilities are, up to isometry,
     (i) $\C\PP^3$ with the Ziller metric \cite{Zi},
     (ii) $\SU(3)/T^2$,  with the Killing form metric or the invariant
     K\"ahler-Einstein metric(s), (iii) some left-invariant Einstein metric on
     $\SU(2) \times \SU(2)$.

     In case (iii), \cite{NiR} shows further that the only left-invariant Einstein metrics
      which are {\em also} $\Ad(S^1)$ invariant for some circle in $\SU(2) \times \SU(2)$
      must be isometric to the product metric or the Killing form metric on
      $(S^3 \times S^3 \times S^3)/\Delta S^3 \approx \SU(2) \times \SU(2).$ It is quite
      likely these are the only left-invariant Einstein metrics, up to isometry. It is
      interesting also to note that as a homogeneous space,  $\SU(2) \times \SU(2)$
      is not of finite type in B\"ohm's sense. Its simplicial complex consists of two
      points.

{\bf (b) Seven-dimensions:} (\cite{Ni2}) The assumptions in this classification are the
     same as those in the $6$-dimensional case. Besides symmetric spaces, the possibilities
     are, up to isometry,

    (i) $\Sp(2)/\Sp(1)$ where the subgroup is embedded by the $4$-dimensional symplectic
    complex representation of $\Sp(1)$ (this space is isotropy irreducible),

    (ii) the Stiefel manifold $\SO(5)/\SO(3)$ with the Kobayashi metric \cite{Ko},

    (iii) $S^7 = \Sp(2)/\Sp(1)$ with the Jensen metric \cite{Jen2},

    (iv) the Aloff-Wallach spaces $\SU(3)/\U_{pq}$ with $p, q$ relatively prime and one of
    two non-isometric homogeneous Einstein metrics \cite{Wa1}, \cite{PP1},

    (v) all simply connected circle bundles over $\C\PP^2 \times \C\PP^1$, each with
      a unique $\SU(3) \times \SU(2)$ invariant metric \cite{CaDF}, \cite{WaZ2},

    (vi) all simply connected circle bundles over $\C\PP^1 \times \C\PP^1 \times \C\PP^1$,
      each with a unique $\SU(2) \times \SU(2) \times \SU(2)$-invariant metric that
      is non-symmetric \cite{DFV}, \cite{Ro}, \cite{WaZ2}.

 An important part of the classification in \cite{Ni2} is the analysis of the Einstein
 metrics on the two exceptional members of (iv), cf Examples \ref{Wallach1} and  \ref{Wallach2}.

 It should be noted that a subfamily of the Einstein manifolds in (vi) consists of the
 product with $\C\PP^1$ of $5$-dimensional homogeneous Einstein spaces \cite{AlDF}. These
 consist of all simply connected circle bundles over $\C\PP^1 \times \C\PP^1$, an exceptional
 member of which is $\SO(4)/SO(2)$ (the unit tangent bundle of $S^3$), cf Example \ref{unitang}.

 In terms of bundle constructions, (v) and (vi) above belong to a large
 family of (mostly non-Sasakian) bundle type Einstein metrics on torus bundles over
 a product of Fano K\"ahler-Einstein manifolds found in \cite{WaZ2}.

\smallskip

{\bf (c) Homogeneous spaces with two irreducible isotropy summands:} (\cite{DiK})
    This classification is done under the assumptions that $G$ is a compact connected
    simple Lie group, $K$ is a closed subgroup,  $G/K$ is simply connected, and the
    isotropy representation consists of exactly two real irreducible summands.
    It follows that the graph and simplicial complex of $G/K$ are not particularly useful
    tools in this case. The space of unit volume invariant metrics
    is $1$-dimensional except when the two irreducible summands
    are equivalent, i.e., when $G/K = \Spin(8)/G_2 \approx S^7 \times S^7$. The invariant Einstein
    metrics on this exceptional space have been determined in \cite{Ke}. Up to isometry,
    these are the product metric and the Killing form metric.

    For the rest of the cases, the Einstein condition reduces to a cubic polynomial
    equation (in one variable), which becomes a quadratic equation when there is a
    non-trivial intermediate subalgebra between $\kf$ and $\g$. When $\kf$ is a maximal
    subalgebra, there is always an Einstein metric \cite{WaZ1}, but it has not been
    determined exactly when there would be $3$ Einstein metrics. In the non-maximal case,
    existence is decided by the discriminant of the quadratic. This is
    examined for all cases in \cite{DiK}, although some cases were treated earlier
    in \cite{WaZ1}. There seems to be no obvious pattern for existence/non-existence,
    despite the large number of spaces involved.

 \smallskip

{\bf (d) Below dimension $13$:} In (\cite{BoK}) the existence of homogeneous
     Einstein metrics is examined for every compact simply connected homogeneous space
     of dimension less than $13$. From this one concludes that each such homogeneous
     space with dimension less than $12$ admits at least one homogeneous
     Einstein metric. However, it leaves open the problem of finding {\em all}
     homogeneous Einstein metrics on some of the homogeneous manifolds, e.g., the Lie groups
     $\Spin(4), \SU(3), \Sp(2)$, $\Sp(2)/\U_{p,q}$, and $\SU(4)/S(\U(2)\cdot \U_{p,q})$
     where $\U_{p,q}$ denotes a circle lying in a $2$-torus.

     In dimension $12$ it is known from \cite{WaZ1} that $\SU(4)/\Sp(1)$
     has no homogeneous Einstein metrics, where $\Sp(1)$ is embedded into $\SU(4)$ via its
     symplectic complex representation of dimension $4$. Two new families with
     no $G$-invariant Einstein metrics are found in \cite{BoK}: (i) $K=\Delta \SU(2)\cdot \U_{p,q} \subset
     {\rm S}(\U(2) \U(1))\times {\rm S}(\U(2) \U(1)) \subset \SU(3) \times \SU(3)=G$ with $p\neq q$,
     and (ii) $K = \Delta \Sp(1) \times \Sp(1) \subset \SU(2) \times (\Sp(1) \times \Sp(1))
     \subset \SU(3) \times \Sp(2)=G$.

     Some of the new Einstein metrics found in \cite{BoK} require the use of the
     general existence theorems from \cite{BWZ} and \cite{Bo1}. Others still require
     studying the Einstein equations directly with the help of computers, e.g.,
     $\SU(4)/\U(2)$, where $\U(2) \subset \SO(4) \cap \Sp(2) \subset \SU(4)$, cf \cite{Ni1},
     $(\SU(3) \times \SU(3))/ (\Delta\Sp(1) \U(1))$ where $\Delta\Sp(1) \subset \SO(3) \times \SU(2)$,
     and $(\SU(3) \times \Sp(2))/(\Delta\Sp(1)\times \Sp(1))$ where
        $\Delta\Sp(1)\times \Sp(1) \subset (\SO(3) \times \Sp(1))\times \Sp(1)$.

\smallskip

Besides classification theorems, attempts have also been made to show that certain
familiar homogeneous manifolds $G/K$ admit many invariant Einstein metrics. A common
strategy is to look among $G$-invariant metrics which have additional symmetries, e.g.,
invariance with respect to the adjoint action of some closed subgroup of the gauge
group $N_G(K)/K$. This strategy simplifies the Einstein condition by restricting
attention to a smaller family of invariant metrics. For example, Einstein left invariant
metrics on compact Lie groups were produced by \cite{Jen2} and \cite{DZ} in this way.
We mention here some recent results pertaining to coadjoint orbits, Stiefel manifolds,
and homogeneous fibrations.

\smallskip

{\bf (1) Coadjoint orbits}: These are also called generalized flag manifolds
or K\"ahlerian $C$-spaces in the literature. For a compact connected semisimple
Lie group $G$, the adjoint representation is self-contragredient, so there is
no distinction between orbits in the adjoint and the coadjoint representations.
Each such orbit is of the form $G/C(T)$ where $C(T)$ denotes the centralizer of
a torus in $G$. The invariant almost complex structures are all integrable and
have positive first Chern class. It is classical \cite{Be} that each invariant
complex structure has an associated K\"ahler-Einstein metric which is unique up
to complex automorphisms. Note that $G$ and $C(T)$ have the same rank, so the
isotropy representation always splits into pairwise inequivalent irreducible summands.
Hence they are prime homogeneous spaces in B\"ohm's sense.

There has been a lot of work trying to find non-K\"ahler invariant Einstein metrics
on coadjoint orbits, and even to determine the set of all invariant Einstein metrics.
The article \cite{ArCh2} gives an up-to-date account of this activity.
For all coadjoint orbits that have been examined, non-K\"ahler invariant Einstein
metrics were found to exist.

One theme is to look at coadjoint orbits whose isotropy representation splits into
a small number of irreducible summands. The case of $2$ irreducible summands
was analysed in \cite{ArCh1}, while the case of $3$ irreducible summands was
examined in \cite{Ar} and \cite{Ki}. In these works, all invariant Einstein metrics
were determined. The paper \cite{ArCh2} treats the case of $4$ irreducible summands
as well as isometry issues. All the invariant Einstein metrics were determined for
these spaces up to isometry except for the two infinite families
$\SO(2m)/(\U(p) \times U(m-p))$  with $m \geq 4, \, 2 \leq p \leq m-2$
and $\Sp(m)/(\U(p) \times U(m-p))$ with $m \geq 2, \, 1 \leq p \leq m-1$, for which
some new Einstein metrics were also found. The complete determination of invariant
Einstein metrics was achieved in \cite{ArChS1} for the first family. Except for
ten members of this family (with small $m$ and $p$), there are exactly two non-K\"ahler
Einstein metrics. For the second family, see \cite{ArChS2} where it is shown that
there are exactly two non-K\"ahler Einstein metrics.

A second theme is to study the full flag manifolds, i.e., the principal coadjoint
orbits. If $G$ is a simply-laced compact connected simple Lie group, then the Killing
form metric is Einstein, as was observed by Ziller and myself in the distant past.
For $\SU(n)/T$ (where $T$ is a maximal torus), further non-K\"ahler Einstein metrics
were first found in \cite{Ar}. Other solutions to the Einstein equations were found
by Sakane \cite{Sa}, and most recently by Dos Santos and Negreiros \cite{DoNe}.

\smallskip

{\bf (2) Stiefel manifolds}: The real Stiefel manifolds $\SO(m+pn)/SO(m)$ with
$p > 1, m>n\geq 3$ are shown in \cite{ArDNi2}, \cite{ArDNi3} to have at least
four $\SO(m+pn) \times \SO(n)^p$-invariant Einstein metrics, two of which were
constructed first in \cite{Jen2} and the rest are new. It is also shown in these
works that given any positive integer $q$ one can find a real Stiefel manifold with
at least $q$ invariant Einstein metrics. Analogous results are proved in \cite{ArDNi1}
for the quaternionic Stiefel manifolds $\Sp(m+pn)/\Sp(m)$ with $p>1, m \geq n\geq 1$.

\medskip

We conclude this section with a discussion of compact homogeneous spaces without
any invariant Einstein metrics. Of course it is entirely possible that the
underlying smooth manifolds admit Einstein metrics with less or even no symmetries.
In dimensions $\geq 5$, it is still unknown whether every closed manifold admits an
Einstein metric or not, and there are no known general obstructions to existence.
For positive Einstein metrics, the fundamental group of the manifold must be finite,
by Bonnet-Myers, and in the spin case, all obstructions to positive scalar curvature
must also vanish.

Coming back to the homogeneous case, in view of the discussion at the beginning
of this section, we may assume that $G/K$ is simply connected with $G$ compact,
simply connected, and semisimple. The first non-existence examples start in dimension
$12$ (cf \S 3(d) above). An immediate consequence of Theorems \ref{graphthm} and
\ref{simcplxthm} is that for $G/K$ to have no $G$-invariant Einstein metrics,
its simplicial complex must be contractible and its graph must be connected.

Sufficient conditions for non-existence were studied in \cite{Bo2}. These are
based on the following observation. Consider the decomposition of $\p$ given in
(\ref{decomp}), where we assume there are at least two summands. For a $G$-invariant
metric $g$, its traceless Ricci operator $\hat{r} - \frac{S}{n}q$
(cf Remark \ref{recover}) preserves this
decomposition by invariance considerations. Since this operator is zero when $g$
is Einstein, it follows that if for all $G$-invariant metrics and some summand $\p_i$,
the component of the traceless Ricci tensor in this summand is always positive
definite (resp. negative definite), then $G/K$ cannot admit a $G$-invariant
Einstein metric.

B\"ohm then analyses the structure of homogeneous spaces which satisfy the
above two definiteness criteria and obtains

\begin{thm} $($\cite{Bo2}, Theorems B, C$)$ \label{structure}
Let $G/K$ be a compact homogeneous space with finite fundamental group.
Suppose all $G$-invariant metrics on $G/K$ have the property that
for a fixed summand $\p_i$ in $($\ref{decomp}$)$, the $\p_i$-part of the
the traceless Ricci tensor is always definite.

$($i$)$ In the positive definite case, $K$ is contained in some compact
subgroup $H$ such that $H/K$ is isotropy irreducible of dimension $\geq 2$,
and all $G$-invariant metrics on $G/K$ are of Riemannian submersion type with
respect to the fibration $H/K \rightarrow G/K \rightarrow G/H$.

$($ii$)$ In the negative definite case, $K$ is contained in some compact
subgroup $H$ such that $G/H$ is isotropy irreducible of dimension $\geq 2$,
and $(\h, \kf) = \oplus (\h_j, \kf_j)$ where each pair $(\h_j, {\kf}_j)$
corresponds to a compact connected isotropy irreducible space.
\end{thm}

Examples of non-existence are then obtained in \cite{Bo2} by explicitly
examining the Ricci tensors of invariant metrics on coset spaces
$G/K$ of the type described in Theorem \ref{structure} above. Furthermore,
it is possible to build large families of new non-existence examples out of
known ones by a gluing procedure (cf Theorem 4.7 in \cite{Bo2}). In particular,
one obtains in this way simply connected non-existence examples whose spaces
of invariant metrics have arbitrarily large dimension.

\section{\bf Cohomogeneity One Einstein Equations}

A smooth manifold $M^{n+1}$ is said to be of {\em cohomogeneity one} if there is
a Lie group $G$ acting (properly) on it such that the principal orbits are hypersurfaces,
i.e., the orbit space is $1$-dimensional. In the context of Riemannian geometry,
$G$ is a group of isometries of some metric $\overline{g}$ on $M$, and it is usually
assumed to be compact. On the other hand, in complex geometry, it is natural to
take $G$ to be a complex Lie group acting with an open orbit, and we arrive at
the concept of an {\em almost homogeneous complex manifold}, cf \cite{HS}.

Theoretical physicists have of course studied cohomogeneity one Lorentz metrics
ever since the birth of General Relativity. Indeed the non-K\"ahler Page
metric on $\C\PP^2 \sharp \overline{\C\PP^2}$, the first known inhomogeneous positive
Einstein metric on a closed manifold, was produced via a Wick rotation from
the Taub-NUT metric in Relativity \cite{P}. The connection with spatially homogeneous
Lorentz metrics adds a further dimension to the study of cohomogeneity one Einstein
manifolds that is absent in the study of the cohomogeneity one condition for other
types of geometries.

The mathematical formulation of the Einstein condition for cohomogeneity one metrics
was first undertaken in \cite{BB}. The role of the second Bianchi identity was
pointed out in \cite{Ba}. The initial value problem at a singular orbit was
posed and solved in many cases in \cite{EW}. The Hamiltonian formulation of the
Einstein equations was given in \cite{DW3}. (See appendix E in \cite{Wal}
for the Hamiltonian formulation of General Relativity.) We will now describe
these developments in this section. In the next section we will address integrability
issues and the search for interesting subsystems of the cohomogeneity one Einstein
equations.

\smallskip

\noindent{\bf 4A. Basic set-up}

\smallskip

Throughout \S 4 - \S 5, we will let $G$ be a compact Lie group acting by isometries
on a Riemannian manifold $(M, \overline{g})$ of dimension $n+1 \geq 4$ with $1$-dimensional
orbit space $I$. Let us choose a base point $x_0 \in M$ lying on a principal orbit $P$, and
denote by $K$ its isotropy group. We next choose a unit speed geodesic $\gamma$ which
passes through $x_0$ and intersects all principal orbits orthogonally. The (arclength) parameter
along $\gamma$ will be identified with the coordinate in $I$. Since we are interested in
the Einstein condition, we will assume that the orbit space of the $G$-action is not a circle.
(Otherwise, $M$ would be compact with infinite fundamental group and cannot support an
invariant non-flat Einstein metric.)

The choice of $\gamma$ provides us with a diffeomorphism
$$\phi: I^0 \times P \longrightarrow M_0$$
where $M_0$ denotes the open and dense set of points in $M$ lying on  principal
orbits and $I^0$ denotes the interior of the orbit space $I$. Explicitly,
$\phi(t, a[K]) = a \cdot \gamma(t)$ and
\begin{equation} \label{metricform}
\phi^* \overline{g} = dt^2 + g_t
\end{equation}
where $g_t$ is regarded as a curve in ${\mathscr M}^G(G/K)$ in the notation of \S 1.
By fixing a background $G$-invariant metric $b$ as in \S 1, we obtain via (\ref{gtoq})
a corresponding path $q_t$ in ${\mathscr C} := S_{+}(\p)^K$, the space of
$b$-symmetric, positive definite, $\Ad_K$-invariant automorphisms of $\p$.
We also let $\mathscr S$ denote $S(\p)^K$, the $b$-symmetric $\Ad_K$-invariant
endomorphisms of $\p$. We will regard ${\mathscr C} \subset {\mathscr S}$
as the {\em configuration space} in our Lagrangian formulation of the cohomogeneity one
Einstein equations. The {\em velocity phase space} is then $T{\mathscr C} \approx
{\mathscr C} \times {\mathscr S}$, which has coordinates $(q, \dot{q})$ where $\dot{q}$
is given by
$$ \dot{g}(X, Y) = b(\dot{q}(X), Y).$$

In order to write down the Einstein equations, we introduce the unit normal field
$N:=\frac{\partial}{\partial t}$ and the {\em shape operator} $L$ of the principal orbits
given by
$$ L(X) = \overline{\nabla}_X N $$
for $X$ tangent to the principal orbits. ($\overline{\nabla}$ is the Levi-Civita connection
of $\overline{g}$.) We regard $L_t$ (for fixed $t$) as a $G$-equivariant endomorphism of
$TP$ that is symmetric with respect to $g_t$ (but not necessarily with respect to $b$).
By means of the Gauss and Codazzi equations (cf \cite{EW}, \S 2) we obtain the Einstein system
\begin{eqnarray} \label{co1eq}
\dot{q} & = & 2 q \circ L,  \label{eq1} \\
-\dot{L} - ({\rm tr} L) L + \varepsilon r   & = & \varepsilon \Lambda \I,  \label{eq2} \\
-{\rm tr}(\dot{L}) - {\rm tr}(L^2)   & = &  \varepsilon \Lambda, \label{eq3} \\
d({\rm tr} L) + \delta^{\nabla} L & = & 0 \label{eq4}
\end{eqnarray}
where $r$ is the Ricci endomorphism of $g$ defined by ${\rm Ric}(g)(X, Y) = g(r(X), Y),$
$\I$ is the identity endomorphism of $T(G/K)$, and $\delta^{\nabla} : \Omega^1(G/K, T(G/K))
\longrightarrow T(G/K) \approx T^*(G/K)$ is the divergence operator followed by
the isomorphism induced by the metric $g$. The constant $\varepsilon$ equals $1$
in the Riemannian situation that we are considering. To obtain the Einstein condition for
a Lorentz metric of cohomogeneity one of the form $-dt^2 + g_t$, one only needs to
set $\varepsilon = -1$ in the above.

Equation (\ref{eq1}) may be viewed as defining $L$. Equation (\ref{eq2}) represents
the Einstein condition along the principal orbits, while equation (\ref{eq3}) represents
the Einstein condition in the direction of $N$. Equation (\ref{eq4}) expresses the
Einstein condition in the mixed directions. If we take the trace of Equation (\ref{eq2}) and
subtract from it Equation (\ref{eq3}), we obtain the {\em conservation law}
\begin{equation} \label{conlaw}
\tilde{\sf H}:=  (\tr L)^2 - \tr(L^2)  - \varepsilon \,S +(n-1)\,\varepsilon \Lambda = 0.
\end{equation}
We shall see that $\tilde{\sf H}$ is closely related to the Hamiltonian that arises
in the formulation of equations (\ref{eq1})-(\ref{eq4}) as a Hamiltonian system
(cf (\ref{ham})).

Note that Equation (\ref{eq4}) together with $\tilde{\sf H} =0$ constitute the
analogue of the Einstein constraint equations in General Relativity.

We now focus on the Riemannian case and take $\epsilon=1$ in the following.
We want to interpret the Einstein equations (\ref{eq1})-(\ref{eq4}) as a dynamical
system in the configuration space $\mathscr C$. To do so, we first replace (\ref{eq3})
by the constraint (\ref{conlaw}). To handle equation (\ref{eq4}), we need to
bring in singular orbits.

We will assume that our cohomogeneity one manifold $M$ has at least one singular orbit,
i.e., the orbit space $I$ is either a half-open interval or a closed interval.
In the former case, we have one singular orbit, $M$ is non-compact, and we are interested
in complete Einstein metrics on $M$. In the latter case, there are two singular orbits
and $M$ is a closed manifold. Without loss of generality we may assume that $t=0$ is
an endpoint of the orbit space $I$. Let $Q$ denote the singular orbit at $t=0$. Denote
by $H$ the isotropy group of $\gamma(0)$, so that $Q=G/H$. It follows that we have
$K \subset H \subset G$ and $H/K \approx S^k$. As is well-known, a $G$-invariant tubular
neighbourhood of $Q$ has the form $G \times_{H} D^{k+1}$ where the linear action of $H$
on the disc $D^{k+1}$ is the slice representation at $\gamma(0)$, and every normal sphere
subbundle of $G \times_{H} D^{k+1}$ is a principal orbit.

A smooth solution of Equations (\ref{eq1})-(\ref{eq4}) in general gives an Einstein
metric only in the open subset $M_0 \subset M$ via (\ref{metricform}). Additional
boundary conditions have to be satisfied in order for this metric to extend over the
singular orbits to a smooth metric on $M$. These {\em smoothness conditions} have been
completely worked out in \S 1 of \cite{EW}. They involve studying tensor invariants
of the slice and isotropy representations at the singular orbits. It should be
remarked that these conditions can be quite complicated, cf Example 2.1 in
\cite{Wa3} or \cite{BaHs}, and are different from the types of boundary conditions
usually studied in the ODE literature.

\begin{prop} $($A. Back \cite{Ba} $)$ \label{bianchi}
Let $\overline{g}$ be a $C^3$ Riemannian metric of cohomogeneity one in a
tubular neighbourhood $E:=G \times_{H} D^{k+1}$ of a singular orbit $Q$ as in
the above discussion. Assume in addition that $\dim Q < \dim P$.  If
$($\ref{eq1}$)$ and $($\ref{eq2}$)$ hold on $E \setminus Q$, then
$($\ref{eq3}$)$ and $($\ref{eq4}$)$ also hold, i.e., $\overline{g}$ is Einstein.
\end{prop}

Hence in the situation where at least one singular orbit has dimension strictly
smaller than that of the principal orbits, the Einstein condition becomes a constrained
dynamical system in $\mathscr C$, provided the smoothness conditions at the singular
orbits can be verified.

When one studies Equations (\ref{eq1})-(\ref{eq4}) for specific principal orbit types,
a frequent assumption is for the isotropy representation of $G/K$ to have
no multipliticies. In this situation, an argument of B\'erard Bergery (cf (3.18)
in \cite{BB}) shows that (\ref{eq4}) always holds, even if there is no singular
orbit of strictly smaller dimension. When the isotropy representation of $G/K$
has multiplicities, one can often enlarge $G$ and $K$ by adding finite or toral
groups so that the resulting coset space has multiplicity free isotropy representation.
This amounts to restricting $g_t$ to lie in a subset of the set of invariant
metrics of $G/K$. However, the no multiplicities hypothesis, even in this form,
remains rather restrictive.

When one considers $G/K$ with multiplicities in its isotropy representation, a
technical but important question arises. When does the Einstein condition
force $g_t$ to be simultaneously diagonalizable (for all $t$) with respect to some fixed
($\Ad_K$-invariant) decomposition of $\p$ ? A well-known example from General
Relativity where this happens is the case of the Bianchi IX metrics, i.e., when
$G/K$ is $\SU(2)$. In this case, simultaneous diagonalizability actually follows just from
Equation (\ref{eq4}). A first attempt at analysing the relation between simultaneous
diagonalizability and Equation (\ref{eq4}) has been made in \cite{Da}, where in
particular the principal orbits $\SO(m+2)/SO(m)$ and $\SU(m+2)/\U(m)$ were studied.

\smallskip

\noindent{\bf 4B. Hamiltonian formulation}

\smallskip

We next consider the reformulation of the above dynamical system in Hamiltonian
terms. The starting point is the Hilbert action in the case of a compact $M$.
Recall that the scalar curvature of $\overline{g}$ is given by
$$ \overline{R}= -2\, \tr(\dot{L}) -\tr(L^2) -(\tr L)^2 + S, $$
where $S$ is the scalar curvature of the principal orbits. If we assume
the background metric $b$ on $P=G/K$ has volume 1, then we have
\begin{eqnarray*}
{\mathcal A}(\overline{g}) &=& \int_M \, \overline{R}\, d{\rm vol}_{\overline{g}} \\
     &=& \int_0^{\tau} \, \left(S-\tr(L^2)-(\tr \, L)^2 -2\tr(\dot{L}) \right) v dt
\end{eqnarray*}
where $v$ is the relative volume (cf (\ref{relvol})) and $I=[0,\tau]$. Since
$\dot{v} = (\tr L)v$, upon integration by parts we obtain
$$ {\mathcal A}(\overline{g}) = \int_0^{\tau} \, \left((\tr \, L)^2 - \tr(L^2) + S\right) v dt
        -2 (\tr \,L)v \left|_0^{\tau} \right.. $$
The boundary term vanishes when the singular orbits have codimension $k+1 \geq 3$
since for small $t$ we have $\tr \, L \sim \frac{k}{t}$ and $v(t) \sim t^k$.

The above computation suggests taking the Lagrangian to be
\begin{equation*}
  {\mathsf L} := ((\tr \, L)^2 - \tr(L^2) + S)v -(n-1)\Lambda v
\end{equation*}
where the last term on the right is added because of the volume constraint
in the Hilbert action.

Note that $(\tr \, L)^2 - \tr(L^2) = 2\langle L, L \rangle$ where
\begin{equation} \label{quadform}
   \langle \alpha, \beta \rangle := \frac{1}{2} \left((\tr \,\alpha)(\tr \,\beta) - \tr(\alpha \beta) \right)
\end{equation}
is the natural non-degenerate $\Ad_K$-invariant symmetric bilinear form on ${\rm End}(\p)$.
(Our sign for $\langle \cdot, \cdot \rangle$ is opposite to that in \cite{DW3}.)
It is positive definite on the $b$-skew-symmetric operators of $\p$ and a Lorentz
metric of signature type $(1, n-1)$ on the $b$-symmetric operators of $\p$.
Using (\ref{eq1}) we obtain
\begin{equation} \label{lag}
{\mathsf L}(q, \dot{q}) = \frac{1}{2}\langle q^{-1}\dot{q}, q^{-1}\dot{q} \rangle \,v(q)
       + \left(S(q)-(n-1)\Lambda\right) v(q).
\end{equation}
The first term on the right can be interpreted as kinetic energy and the second
term on the right as potential energy. Thus the kinetic energy depends on the
shape operator and relative volume of the principal orbits, and the potential
energy depends on the Einstein constant, and the scalar curvature function and
relative volume of the principal orbits.

For the Hamiltonian formulation, we let the {\em momentum phase space} to be the
cotangent bundle of $\mathscr C$: $T^*{\mathscr C} \approx {\mathscr C} \times {\mathscr S}^*$,
with coordinates $(q, p)$. The symplectic structure is the canonical one, given by
$$ \omega((q, p), (\tilde{q}, \tilde{p})) = p(\tilde{q}) -\tilde{p}(q).$$
Let $\langle \cdot, \cdot \rangle^*$ denote the symmetric bilinear form on ${\mathscr S}^*$
induced by $\langle \cdot, \cdot \rangle$. Recall also that $\GL(\p)$ acts on $\End(\p)$
on the left, and this induces a dual action on $\End(\p)^*$ which specializes to
$$ (q \cdot p)(\xi):= p(q^{-1}\cdot \xi), \,\,\,\,\, \mbox{\rm for} \,\,\,\,\,
      q \in {\mathscr C}, \, \xi \in {\mathscr S}, \, p \in {\mathscr S}^*. $$

The {\em Legendre transformation} from $T{\mathscr C}$ to $T^*{\mathscr C}$ is given by
\begin{equation} \label{legendre}
  p(\xi) = {\mathsf L}_{\dot{q}} (\xi) = \langle q^{-1}\dot{q}, q^{-1}\xi \rangle \,v.
\end{equation}
It follows that $(q^{-1} \cdot p)(\xi) = v \langle q^{-1} \dot{q}, \xi \rangle$, and so
$$ \langle q^{-1} \cdot p, \, q^{-1} \cdot p \rangle^* = v^2 \langle q^{-1}\dot{q}, \, q^{-1}\dot{q} \rangle
      = 4v^2 \langle L, L \rangle. $$
The corresponding Hamiltonian on $T^*{\mathscr C}$ is therefore
\begin{equation} \label{ham}
{\mathsf H}(q, p) = \frac{1}{2v(q)}\,\langle q^{-1} \cdot p, \, q^{-1} \cdot p \rangle^*
        + \left((n-1)\Lambda -S(q)\right) v(q).
\end{equation}
It will be convenient to denote the quadratic form
$\frac{1}{2}\, \langle q^{-1} \cdot p, \, q^{-1} \cdot p \rangle^*$ by $J(p, p)$.

The relationship between solutions of Hamilton's equation for the Hamiltonian (\ref{ham})
and solutions of (\ref{eq1})-(\ref{eq4}) is then given by
\begin{prop} $($\cite{DW3}, Proposition 1.13 $)$ \label{hamflow}
Let $M$ be a cohomogeneity one $G$-manifold with principal orbit type $G/K$ and a singular
orbit $G/H$ with  $ \dim G/H < \dim G/K$. Suppose $(q(t), p(t))$ is a solution of Hamilton's
equation for the Hamiltonian $($\ref{ham}$)$ on $T^*{\mathscr C}$ that lies on the zero energy
hypersurface $\{\mathsf H = 0\}$. Let $(q(t), \dot{q}(t))$ be the corresponding curve
given by the inverse of the Legendre transformation $($\ref{legendre}$)$. Then $(q(t), \dot{q}(t))$
satisfies $($\ref{eq2}$)$.

If the $G$-invariant metric represented by $(q(t), \dot{q}(t))$ extends over the singular
orbit as a $C^3$ Riemannian metric, then it is actually Einstein.
\end{prop}

\begin{rmk} \label{hypersurface}
It is easily seen that if $G/K$ is not a torus or if $\Lambda \neq 0$, then
$\{\mathsf H = 0\}$ is a regular hypersurface in $T^*{\mathscr C}$. If $G/K$ is
a torus and $\Lambda =0$, then the singular points of $\mathsf H$ consist of
$(q, 0)$ where $q$ is a homogeneous flat metric on $G/K$.
\end{rmk}

\smallskip

\noindent{\bf 4C. Local existence}

\smallskip

Are there any local obstructions to solving the cohomogeneity one Einstein equations ?

Since we are dealing with ordinary differential equations here, there is certainly
no obstruction to solving the equations in a neighbourhood of a principal orbit
provided that the initial metric and shape operator satisfy (\ref{eq4}) and the required
invariance properties. (This is the elementary analogue of the Einstein constraint problem
in General Relativity.) The local solution with prescribed initial data is unique.

Solving the equations in the neighbourhood of a singular orbit is less trivial because
there is a non-regular singular point of (\ref{eq1})-(\ref{eq2}) at each endpoint of
the orbit space $I$. This corresponds to the geometric fact that as the principal orbits
collapse down to the singular orbit, the Ricci tensor and shape operator blow up in the
collapsing (radial) direction. One must also make sure that the solutions to
(\ref{eq1})-(\ref{eq2}) satisfy the additional smoothness conditions for them to
represent smooth metrics near the singular orbit.

\begin{thm} $($\cite{EW}$)$  \label{IVP}
Let $G \times_H D^{k+1}$ be a tubular neighbourhood of the singular orbit $Q=G/H$
in a cohomogeneity one manifold with principal orbit type $G/K$, where $G$ is a compact
Lie group, $K \subset H \subset G$,  and $G/K$ is connected. Assume that as
$K$-representations, the slice representation $D$ and the isotropy representation of $G/H$
do not have any irreducible $\R$-subrepresentations in common. Then given any $G$-invariant
metric on $Q$ and any
$G$-equivariant tensor $\nu(Q) \longrightarrow S^2(T^*Q)$ $($prescribed shape operator$)$,
there exists a $G$-invariant Einstein metric of any prescribed sign $($including $0$$)$
of the Einstein constant on an open subdisc bundle of $G \times_H D^{k+1}$.
\end{thm}

\begin{rmk} \label{pfivp}
In the above $\nu(Q)$ denotes the normal bundle of $Q$ in $M$. Note that $G$-equivariance
means that the prescribed shape operator can be thought of as an $H$-equivariant linear
map $D \longrightarrow S^2(\g/\h)$. It follows that the mean curvature vector is zero,
which is a well-known necessary condition for having a smooth cohomogeneity one metric
(cf \cite{HL}).
\end{rmk}

Theorem \ref{IVP} is proved by constructing a formal power series solution to
(\ref{eq1})-(\ref{eq2}), and then applying Picard iteration to a high order
truncation of the power series. The main difficulty is that the recursion operator
is not surjective. However, it turns out that the right-hand side of the recursion
relation always lies in the range of the recursion operator provided that the smoothness
conditions are imposed in all earlier steps of the recursion. That the recursion operator
is also not injective gives non-uniqueness of the local solution. In fact, examples
given in \cite{EW} show that one cannot bound the level at which the recursion operator
becomes injective (although in any specific instance the operator does become
injective after finitely many steps). In all likelihood, the technical hypothesis
on the slice and isotropy representations at the singular orbit can be removed
by a more careful study of the recursion and the necessary invariant theory.

\begin{rmk} \label{reidegeld}
In a recent paper \cite{Re3}, Reidegeld has determined explicitly all the parameters
that would specify uniquely a local solution to the full cohomogeneity one Einstein equation
when the principal orbit is a generic Aloff-Wallach space and the singular orbit is
one which can occur for metrics with $\Spin(7)$ holonomy. For the exceptional
Aloff-Wallach spaces, he determined the parameters for the subclass of cohomogeneity
one metrics for which the induced metrics on the principal orbit are diagonal
with respect to the decompositions indicated in Examples \ref{Aloff1} and \ref{Aloff2}.
(See \S {\bf 6B} for further details.)
\end{rmk}

\smallskip

\noindent{\bf 4D. Global non-existence}

\smallskip

While there appear to be no local obstructions to solving the cohomogeneity one
Einstein equations, global obstructions do exist, as the following result of
B\"ohm shows.

\begin{thm} $($\cite{Bo3}$)$ \label{nonexist}
Let $M$ be a closed $G$-manifold with cohomogeneity one. Suppose
that $G/K$ is the principal orbit type where $G$ is a compact Lie group
and $G/K$ is connected. Fix a bi-invariant background metric $b$ on $G$.
Let $Q_i= G/H_i,$ $i=1,2$ be the two singular orbits with $K \subset H_i$.
Assume that $\h_i = \kf \oplus \s_i$ are $b$-orthogonal $\Ad_K$-invariant
decompositions, and $\p_0 \oplus \cdots \oplus \p_{\ell}$ is the decomposition
$($\ref{decomp}$)$ of the isotropy representation of $G/K$ described in \S 1.

If for some $j$, $\p_j$ is $\Ad_K$-irreducible, $\p_j \cap (\s_1 \cup \s_2) = \{0\},$
and the restriction of the trace-free part of the Ricci tensor of any
$G$-invariant metric on $G/K$ to $\p_j$ is negative definite, then $M$ does not
admit any smooth $G$-invariant Einstein metric.
\end{thm}

\begin{example} \label{noalmhomog}
Examples of cohomogeneity one manifolds which satisfy the hypotheses of Theorem
\ref{nonexist} include $S^{k+1} \times (G^{\prime}/K^{\prime}) \times M_3 \times \cdots \times M_r$
where $M_3, \cdots, M_r$ are arbitrary compact isotropy irreducible homogeneous spaces,
$G^{\prime}/K^{\prime}= {\rm SU}(\ell + m)/({\rm SO}(\ell) {\rm U}(1) {\rm U}(m))$
(a bundle over a complex Grassmannian with a symmetric space as fibre), and
$\ell \geq 32,\, m=1, 2,\, k=1, \cdots, [\ell/3]$  (cf \cite{Bo3}). The significance
of the spaces $G^{\prime}/K^{\prime}$ is that they do not admit any
$G^{\prime}$-invariant Einstein metrics (cf \cite{WaZ1} and Theorem \ref{structure}(ii))
in \S 3).
\end{example}

\begin{rmk} The cohomogeneity one setting (in both its Lagrangian and Hamiltonian forms)
has also been exploited in the construction of Ricci solitons. Recall that these consist
of a triple $(M, g, X)$, where $X$ is a vector field on the Riemannian manifold $(M,g)$,
satisfying the equation
$$ {\Ric}(g) + \frac{1}{2} {\mathscr L}_X g + \frac{\epsilon}{2} g = 0.$$
(${\mathscr L}$ in the above denotes the Lie derivative and $\epsilon$ is a constant.)
When the vector field $X$ is the gradient of a smooth function $u$, the soliton is called
a {\em gradient} Ricci soliton. Analogues of the results in this section for these structures
can be found in \cite{DW14} and \cite{DHW}.
\end{rmk}

\section{\bf First Order Subsystems, First Integrals, and Painlev{\'e} Analysis}

Since the cohomogeneity one Einstein equations can be formulated as a constrained
Hamiltonian system  (at least in the situation where there is a singular orbit
of strictly smaller dimension), a natural question is whether additional conserved
quantities exist. The presence of such quantities, as well as quantities which
change monotonically along the flow, often help in establishing existence of
solutions. In the context of Liouville integrability, conserved quantities arise
as functions on momentum phase space which Poisson commute with the Hamiltonian.
Because of the constraint condition, we can expect only commutation on a subvariety
of phase space, and so the concept of a first integral must be suitably modified.
One particular situation is when the subvariety is Lagrangian and invariant under
the Hamiltonian flow. One then obtains a first order subsystem of the cohomogeneity
one Einstein equations with one-half of the total degrees of freedom. In this section,
we will describe recent efforts to explore the above issues in the Ricci flat case.
Accordingly we will let the Einstein constant $\Lambda = 0$ throughout this section.

\smallskip

\noindent{\bf 5A. Superpotentials}

\smallskip

We begin with the concept of a {\em superpotential} in the physics literature.
This is a $C^2$ function $u: {\mathscr C} \rightarrow \R$ satisfying the first order
PDE
\begin{equation} \label{superpot}
    {\mathsf H}(q, du_q) = 0.
\end{equation}
It induces in a completely natural way a flow on the configuration space given by
\begin{equation} \label{FOS}
 \dot{q} = 2v(q)^{-1} J^* \nabla u
\end{equation}
which corresponds to a first order subsystem of the Hamiltonian flow on
momentum phase space. To explain the operator $J^*$ that occurs and, more generally,
the relationship between a superpotential and conserved quantities we need to
introduce certain first integrals which are linear in momenta on phase space.

Let $\dim {\mathscr C} = N$. Note that $\mathscr C$ has $N$ global coordinates
since the exponential map is a diffeomorphism from the linear space $\mathscr S$
to $\mathscr C$. After choosing linear coordinates in $\mathscr S$, we regard $p$
as a general point in ${\mathscr S}^*$, represented by a row vector. Let $\alpha$
be a section of $\Aut(T^*{\mathscr C}) \approx {\mathscr C} \times \GL({\mathscr S}^*)$,
and $\beta$ be a section of $T^*{\mathscr C}$. ($\alpha$ acts on the right of
sections of $T^*{\mathscr C}$.) Then
$$ \Phi = p \cdot \alpha + \beta$$
can be viewed as an $\R^N$-valued function on $T^*{\mathscr C}$. We denote by
$\Phi_j, \,1 \leq j \leq N,$ its components, which we regard as candidates for first
integrals. We introduce the zero sets
$$ {\mathscr V}_{\sf H} := \{(q, p): {\sf H}(q, p) = 0 \}$$
and
$$ {\mathscr V}_{\Phi} := \{(q, p): \Phi_1(q, p) = 0, \cdots, \Phi_N(q, p) = 0\}.$$

By Remark \ref{hypersurface}, ${\mathscr V}_{\sf H}$ is a smooth hypersurface
in phase space if $G/K$ is not a torus.

\begin{prop} $($Proposition 1.1, \cite{DW10} $)$ \label{linint}
With notation as above,  there exists a superpotential $u: {\mathscr C} \rightarrow \R$
iff there is a $\R^N$-valued function of the form $\Phi$ on $T^*{\mathscr C}$ such that

$($a$)$ ${\mathscr V}_{\Phi} \subset {\mathscr V}_{\sf H},$

$($b$)$ for all $1 \leq i < j \leq N$, the Poisson brackets $\{\Phi_i, \Phi_j\}$
        vanish on ${\mathscr V}_{\Phi}$.
\end{prop}

\begin{rmk} \label{lagsubmfd}
Properties (a) and (b) in Proposition \ref{linint} imply that the Poisson brackets
$\{\Phi_i, {\mathsf H}\}$, $1 \leq i \leq N$, vanish on ${\mathscr V}_{\Phi}$.
So $\Phi_i$ are ``on-shell" first integrals (which are linear in momenta)
and the Hamiltonian flow is tangent to ${\mathscr V}_{\Phi}$. Note that $\Phi$
defines a section $p=-\beta\alpha^{-1}$ of $T^*{\mathscr C}$ whose graph is
${\mathscr V}_{\Phi}$. Property (b) now says that ${\mathscr V}_{\Phi}$ is a
Lagrangian submanifold of $T^*{\mathscr C}$, parallelized by the Hamiltonian
vector fields $X_{\Phi_i}$. If we pull back the canonical $1$-form $pdq$ via the
section to $\mathscr C$, we obtain a closed $1$-form because of the Lagrangian
condition. This $1$-form is exact since $\mathscr C$ is contractible, thereby
giving us the superpotential $u$.
\end{rmk}

The flow (\ref{FOS}) is just the pull-back of the Hamiltonian flow via the
section $(q, du_q)$ to the configuration space $\mathscr C$: we have
$\dot{q} = {\mathsf H}_p = \frac{2}{v} J^* p^T = \frac{2}{v}J^* \nabla u$,
where $J^*$ is the linear operator associated with the quadratic form $J$
via the linear coordinates we chose for $\mathscr S$, and $\nabla u$ is the
corresponding Euclidean gradient of $u$.

The physicists' notion of a superpotential therefore has a symplectic
interpretation in terms of conserved quantities, and leads naturally to
a first order subsystem of the cohomogeneity one Einstein equations.
This explains the $C^2$ assumption on superpotentials (so that the
induced vector field on $\mathscr C$ is $C^1$ rather than just continuous)
as well as the insistence that the domain of the superpotential is all
of $\mathscr C$ (we do not want to put an arbitrary limit on where the
trajectories are defined).

Now these assumptions fly in the face of conventional wisdom, for the
superpotential equation (\ref{superpot}) is a (time-independent) Hamilton-Jacobi
equation, and, as a nonlinear implicitly defined first order PDE, global
classical solutions are extremely rare. It turns out that in the Ricci-flat case
physicists were able to find solutions to the superpotential equation, and the
associated first order systems in such cases represented the condition of special
holonomy \cite{BGGG}, \cite{CGLP1}-\cite{CGLP7}. In \cite{DW10}-\cite{DW12} a
classification was given, under assumptions to be described, of those principal orbit
types $G/K$ for which the superpotential equation admits a solution of {\em scalar curvature}
type. This classification includes the cases considered by physicists together
with a few more cases. A brief description of the classification follows.

We will assume that $\Lambda=0$, i.e., the cohomogeneity one metric $\overline{g}$
is Ricci-flat. The principal orbits $G/K$ are assumed to be connected homogeneous
spaces with $G$ compact Lie, $K \subset G$ closed, and multiplicity free isotropy
representation. We will use the notation and set-up in Remark \ref{diagonal}.
Let $d_i$ be the dimension of the irreducible summand $\m_i$ and $d=(d_1, \cdots,
d_r)$. In this situation, $\mathscr C$ is just $\R_{+}^r$, where $\R_{+}$ is the
set of positive reals. We introduce exponential coordinates $x_i = \exp q_i$,
$1 \leq i \leq r$, so that $q=(q_1, \cdots, q_r)$ are coordinates for
$\mathscr S \approx \R^r$.

The scalar curvature function for $G/K$ then has the form
\begin{equation} \label{SCF}
    S = \sum_{w \in {\mathsf W}} A_w \, e^{w \cdot q},
\end{equation}
where $\mathsf W$ is a finite subset of $\Z^r \subset \R^r$ (the {\em weight vectors})
depending only on $G/K$ and $A_w$ are nonzero constants (cf (\ref{diagscal})).
We say that a function is of {\em scalar curvature type} if it is a finite
linear combination of exponentials in the $q_i$. We will assume that
the superpotentials $u$ are of scalar curvature type, and write them in the form
\begin{equation} \label{potform}
     u = \sum_{{c} \in {\mathsf C}} \,F_{c} \,\, e^{c\cdot q}
\end{equation}
where $\mathsf C$ is a finite subset of $\R^r$, and $F_{c}$ are nonzero
real constants.

The classification problem involves determining all principal orbit types
$G/K$ with multiplicity free isotropy representation such that there is
a superpotential of type (\ref{potform}). The focus on superpotentials of
scalar curvature type is certainly restrictive from an analytic point of view.
However, it ensures that the superpotentials are globally defined and $C^2$,
and it allows us to use the geometry of convex polytopes in the search for
solutions.

Note that in our setup of the classification problem, the convex hulls of
$\mathsf W$ and $\mathsf C$ are both unknowns. The superpotential equation
becomes the equations, one for each $\xi \in \R^r$,
\begin{equation} \label{convexeqn}
\sum_{a+c = \xi} J(a,c) \ F_{a} F_{c} = \left\{
  \begin{array}{ll}
    A_w  &   \mbox{if} \ \xi = d + w \  \mbox{for some} \ w \in \mathsf{W} \\
     0    &  \mbox{if} \ \xi \notin d + \mathsf W
   \end{array} \right.
\end{equation}
where $a, c \in {\mathsf C}$ and
\begin{equation} \label{Jform}
   J(p,p) = \frac{1}{n-1} \left(\sum_{i=1}^r p_i \right)^2 - \sum_{i=1}^r \,\frac{p_i^2}{d_i}
\end{equation}
is the quadratic form of signature $(1, r-1)$ occurring in (\ref{ham}).

We will assume that $r > 1$; otherwise $G/K$ is isotropy irreducible, in which
case there is always a superpotential but there cannot be any singular orbits. We also
make the genericity assumption that $\dim \,{\rm conv}({\mathsf W}) = r-1$.
This always holds when $G$ is semisimple (see Theorem 3.11 \cite{DW3}). The equations
(\ref{convexeqn}) imply that
$$ \conv \left(\frac{1}{2}(d+{\mathsf W})\right) \subset {\conv}(\mathsf C),$$
where the inclusion is strict iff there are vertices in $\mathsf C$
which are null with respect to $J$.

\smallskip

The classification in the case where the two convex hulls are the same is given by

\begin{thm} $($\cite{DW10}$)$ \label{nonnullthm}
Let $G$ be a compact connected Lie group and $K$ a closed connected subgroup
such that the isotropy representation of $G/K$ decomposes into a sum of $r$
pairwise inequivalent irreducible real summands. Assume that ${\rm conv}(\mathcal W)$
has dimension $r-1$.

Suppose that the cohomogeneity one Ricci-flat equations with principal orbit
$G/K$ admit a superpotential of scalar curvature type $($\ref{potform}$)$ such that all
elements of $\mathsf C$ are non-null. Then the possibilities, up to permutations of
the irreducible summands, are given by
\begin{enumerate}
\item ${\mathsf W} = \{ (-1) \}$ and $G/K$ is isotropy irreducible,

\item ${\mathsf W} = \{(1, -2), (0, -1)\},$  and
       $G/K = (SO(3) \times SO(2))/\Delta SO(2) \approx \SO(3),$

\item $\mathsf{W} = \{ (1, -2), (-1,0), (0, -1) \},$
       and $G/K$ is one of $(SU(3) \times SO(3))/ \Delta SO(3) \approx SU(3)$,
       $(Sp(2) \times Sp(1))/(Sp(1) \times \Delta Sp(1)) \approx S^7$,
$(SU(3) \times SU(2))/(U(1) \cdot \Delta SU(2))\approx SU(3)/U(1)_{11},$
 or $Sp(2)/(U(1)\times Sp(1))= SO(5)/U(2)\approx \C \PP^3, $

\item $\mathsf{W} = \{ (1,-2,0), (1,0,-2), (0,1,-2), (0,-1,0), (0,0,-1) \},$
     and $G/K$ is $S^7$ written in the form
    $(Sp(2) \times U(1))/(Sp(1)\Delta U(1))$,

\item $\mathsf{W}= \{(1,-1,-1), (-1,1,-1), (-1,-1,1), (-1,0,0), (0,-1,0), (0,0,-1)\}$,
      $d=(2,2,2)$  and $G/K = SU(3)/T$,

\item $ \mathsf{W} = \{ (1,-2,0,0), (1,0,-2,0),
 (1,0,0,-2),(0,1,-1,-1), (0,-1,1,-1),(0,-1,-1,1),$ \\
 $ (0,-1,0,0), (0,0,-1,0), (0,0,0,-1) \},$ and $G/K$
 is an Aloff-Wallach space $SU(3)/U_{pq}$, where
 $U_{pq}$ denotes the circle subgroup consisting of the diagonal
 matrices ${\rm diag}(e^{ip\theta}, e^{iq\theta}, e^{im\theta})$
 with $p+q+m = 0, (p,q)=1$, and $\{p, q, m\} \neq \{1, 1, -2\}$
or $\{1, -1, 0\}$,

\item a local product of an example in $($1$)$ $(n>1)$, $($3$)$, or $($5$)$
     with a circle.
\end{enumerate}
In all of the above cases, there is a superpotential of scalar curvature
type that is unique up to an overall minus sign and an additive constant.
\end{thm}

\begin{rmk} \label{rmk-nonnull}

(a) In Theorem \ref{nonnullthm}, the first order subsystem resulting from case (2)
corresponds to the hyperk\"ahler condition.

(b) The first order subsystems of the second and fourth subcases in (3) as well
as those of (4), (5), (6) correspond to special holonomy $G_2$ or ${\rm Spin}(7)$
according to whether $n=6$ or $7$. These cases will be discussed in further
detail in \S 6.

(c) The first and third subcases of (3) (cf Example \ref{Wallach2}) do not allow the
addition of a singular orbit and are not related to special holonomy.

(d) Case (7) results from a general property of superpotentials of scalar
curvature type without null weight vectors associated with a principal orbit
$G/K$ having no trivial summands in its isotropy representation (cf Remark 2.8
in \cite{DW10}).
\end{rmk}

When there are null vertices in $\mathsf C$ we have

\begin{thm} $($\cite{DW11}, \cite{DW13}$)$ \label{nullthm}
Let $G$ be a compact connected Lie group and $K$ a closed connected
subgroup such that the isotropy representation of $G/K$ is the direct
sum of $r$ pairwise inequivalent $\R$-irreducible summands.
Assume that  $\dim {\conv}(\mathsf W) = r-1$.

Suppose the cohomogeneity one Ricci-flat equations with $G/K$ as
principal orbit admit a superpotential of scalar curvature type $($\ref{potform}$)$
where $\mathsf C$ contains a null vertex. Then either $r \leq 3$, or, up to permutations
of the irreducible summands, we have
\begin{eqnarray*}
 {\mathsf W} &=& \{ (-1)^i, (1^1, -2^i): 2 \leq i \leq r \}, \ \
   d_1=1, \\
 {\mathsf C} &=& \frac{1}{2}(d+ \{(-1^1), (1^1, -2^i): 2 \leq i \leq r \})
  \ \  \mbox{\rm with} \ \  r \geq 2,
\end{eqnarray*}
and the superpotential of scalar curvature type is unique up to an overall minus
sign and an additive constant.
\end{thm}

In the above theorem, $(-1)^i$ denotes the vector whose only nonzero component
is $-1$ occurring in the $i$th position. Similarly, $(1^1, -2^i)$ denotes the vector
whose only nonzero components are $1$ in position $1$ and $-2$ in position $i$.

\begin{rmk} \label{CYcase}
The possibility described in Theorem \ref{nullthm} is realised by circle bundles
over a product of $r-1$ Fano (homogeneous) K\"ahler-Einstein manifolds.
The corresponding first order subsystem corresponds to the Calabi-Yau condition,
cf \cite{CGLP6} and Example 8.1 in \cite{DW10}. The complete Calabi-Yau metrics
were constructed earlier in \cite{WaW}, \cite{DW1}, and Theorem 3.2 in \cite{Wa3}.
\end{rmk}

For the special cases of $r=2, 3$ with a null vector in $\mathsf C$  we have

\begin{thm} $($\cite{DW12}$)$ \label{smallrthm}
Let $G$ be a compact Lie group and $K$ be a closed subgroup such that
$G/K$ is connected and $\kf$ is not a maximal $\Ad_K$-invariant subalgebra
of $\g$. Assume that the isotropy representation of $G/K$ splits into $2$ or $3$
pairwise distinct irreducible real summands and $\mathsf C$ contains a null vector.

If the Ricci-flat cohomogeneity one Einstein equations with $G/K$ as principal
orbit admit a superpotential of scalar curvature type, then, up to permutations of the
irreducible summands, the possibilities are
\begin{enumerate}
\item ${\mathsf W} = \{(-1, 0, 0), (0, -1, 0), (0, 0, -1)\}$ with
       $ d=(3,3,3), (2, 4, 4),$ or $(2, 3, 6)$.
\item ${\mathsf W} = \{(0, -1, 0), (0, 0, -1), (1, -2, 0), (1, 0, -2) \}$ with
      $d_1 =1$,
\item ${\mathsf W } = \{(-1, 0), (0, -1)\}$, with $\frac{4}{d_1} + \frac{1}{d_2} = 1$,
\item ${\mathsf W } = \{(0, -1), (1, -2)\}$, with either $d_1 = 1$ or $\frac{4}{d_1}+
     \frac{9}{d_2} = 1$.
\end{enumerate}
In each of the above cases, there is a superpotential of scalar curvature type that is
unique up to a sign and an additive constant.
\end{thm}

\begin{rmk} \label{rmksmallr}
(a) Case (2) and the $d_1=1$ subcase of Case (4) in the above theorem are
respectively just the $r=3, 2$ cases of the Calabi-Yau case in Theorem \ref{nullthm}.
Furthermore, the second possibility is realised by the complete, non-compact
B{\'erard} Bergery examples \cite{BB}.

(b) Case (3) and the second subcase of Case (4) are realised by the {\em explicit}
doubly-warped examples studied in \cite{DW2}. See Examples 8.2 and 8.3 there for more
details on the superpotentials. Case (1) is realised by the triply-warped examples
studied in \cite{DW5} and \cite{DW10}. The first order
subsystems for the $d=(3,3,3)$ and $(2,4,4)$ are integrable by quadratures. Further
details can be found in Example 8.4 of \cite{DW10}.
\end{rmk}

\begin{rmk} \label{disconncase}
Note that in Theorems \ref{nonnullthm} and \ref{nullthm} we assume that both
$G$ and $K$ are connected, while this is not assumed in Theorem \ref{smallrthm}.
Additional principal orbit types with superpotentials of scalar curvature type can indeed
occur if we drop the connectedness assumption. Examples include $G/K = {\rm O}(3)/({\rm O}(1)
\times {\rm O}(1) \times {\rm O}(1))$ (cf \cite{CGLP5}) and
$G/K = ([SU(2)\times SU(2)\times \Delta U(1)] \ltimes \Z_2 )/ (\Delta U(1) \times \Z_2)
     \approx S^3 \times S^3$ (cf \cite{BGGG} and \cite{CGLP2}).
\end{rmk}

We end this section by pointing out that superpotentials of scalar curvature type
are also known to exist when $\Lambda$ is nonzero. In this situation it is shown
in \cite{DW10}, \S 10 that there must be null vertices in $\conv({\mathsf C})$. Examples
given there include the B{\'e}rard Bergery and Bianchi IX principal orbit
types (cf \cite{CGLP7}).

\smallskip

\noindent{\bf 5B. Generalized first integrals}

\smallskip

Instead of superpotentials for the Hamiltonian (\ref{ham}) one can ask for
the existence of {\em generalized} first integrals which have a polynomial
dependence on momenta (cf chapter VIII of \cite{Koz}). Let $G/K$ be again
a connected principal orbit such that $G$ is compact, $K$ is closed, and
the isotropy representation has no multiplicities. We now seek a solution to
the equation
\begin{equation} \label{polygenFI}
     \{ {\mathsf F}, {\mathsf H} \} = \phi {\mathsf H}
\end{equation}
where $\mathsf F$ and $\phi$ are both finite sums of the form
\begin{equation} \label{intform}
 \sum_b \, A_b(p) \, \exp(b \cdot q)
\end{equation}
with $b \in \R^r$ and $A_b$ polynomial in $p_1, \cdots, p_r$.

One certainly does not expect generalized first integrals to exist frequently.
In order to state a result to this effect, let us recall the set $\mathsf W \subset \R^r$
that arose in the scalar curvature formula (\ref{SCF}). A vector $w \in {\mathsf W}$
is said to be {\em indecomposable} if the only way $kw$ can be written as a sum of
elements in $\mathsf W$ is $w + \cdots + w$ ($k$ times). Also recall the quadratic
form $J$ (see (\ref{Jform})) which occurred in the Hamiltonian (\ref{ham}).

\begin{thm} $($Theorem 4.23, \cite{DW3}$)$ \label{GFI}
With assumptions as above, suppose further that there are at least $3$ irreducible
summands in the isotropy representation of $G/K$, and there exists an ordered basis
$\{w^{(1)}, \cdots, w^{(r)}\} \subset {\mathsf W}$ consisting of indecomposable
vectors such that for each $i$
$$ J(d+w^{(i)}, d+w^{(i+1)}) \neq 0$$
and $\{d+w^{(i)}, d+w^{(i+1)}\}$ spans a time-like $2$-plane in $(\R^r,J).$
Then the only polynomial generalized first integrals are, up to an additive constant,
elements of the ideal generated by the Hamiltonian.
\end{thm}

\begin{rmk} \label{ex-GFI}
 A sufficient condition for the existence of a basis of indecomposable vectors
in $\mathsf W$ is that $G$ is semisimple. One can then check that when $r \geq 3$
the remaining hypothesis in Theorem \ref{GFI} is satisfied provided that
each irreducible summand in the isotropy representation has dimension $\geq 5$.
\end{rmk}

Theorem \ref{GFI} is sharp in the sense that when $r=2$ there are examples with
non-trivial generalized first integrals.

\begin{example} $($\cite{DW2}$)$   \label{Bergery}
Let $G/(K\cdot \U(1))$ be a compact irreducible Hermitian symmetric space
If it is not the hyperquadric $\SO(m+2)/(\SO(m) \SO(2))$, then $G/K$ has
isotropy representation of the form $\p = \I \oplus V$, where the irreducible
summand $V$ can be thought of as the tangent space to the Hermitian symmetric
space at the base point. (For the case of the hyperquadric, see Example \ref{unitang}.)
We have $d_1=1, d_2 = \dim V=n-1$, and the scalar curvature function has the form
$$ S = A_1 \exp(-q_2) - A_2 \exp(q_1 -2q_2)$$
with $A_i > 0$. So ${\mathsf W} = \{(0, -1), (1, -2)\}.$ It is shown in \S 1
of \cite{DW2} that
$$ {\mathsf F} = \frac{p_1^2}{n-1}  \exp(-2q_1+(2-n)q_2) + A_2 \exp(-q_2) $$
is a generalized first integral for the cohomogeneity one Einstein equations
(even when $\Lambda \neq 0$). This gives a conceptual explanation of the
explicit integrability of these Einstein equations, cf \cite{BB}, \cite{PP2}.
\end{example}

Further generalized first integrals (quadratic in momenta) were discovered
in \cite{DW2}.

\begin{example} $($\cite{DW2}$)$  \label{2warp}
Let the principal orbit $P$ be the product of two isotropy irreducible
spaces $G_1/K_1$ and $G_2/K_2$. Then the scalar curvature function takes the
form
$$ S = A_1 \exp(-q_1) + A_2 \exp(-q_2)$$
where $A_i > 0$. So ${\mathsf W} = \{(-1, 0), (0, -1)\}$.
Let $d_i = \dim G_i/K_i$. Then when $(d_1, d_2)=
(2, 8), (3, 6)$ or $(5, 5)$, the cohomogeneity one Ricci flat equations
admit a non-trivial generalized first integral, given respectively by
$$F = -\frac{1}{72} \left( 2p_1-p_2 \right)^{2} e^{-\frac{1}{2}q_1 -7 q_2}
      + A_1 e^{\frac{1}{2} q_1 + q_2},$$
$$F = -\frac{1}{24} (p_1-p_2)^2 e^{-q_1-5q_2} + A_1 e^{q_1+q_2},$$
$$F = -\frac{1}{45} \left( p_1-2p_2 \right)^{2} e^{-2 q_1 -4 q_2} +
      A_1 e^{2q_1+  q_2}.$$
Note that the above pairs $(d_1, d_2)$ are precisely the positive integral
solutions of the diophantine equation $\frac{1}{d_1} + \frac{4}{d_2} = 1$
we encountered in case (3) of Theorem \ref{smallrthm}.

Using the above integrals, one can solve the cohomogeneity one Ricci-flat equations
explicitly to obtain complete Ricci-flat metrics on $(G_1/K_1) \times \R^{d_2 +1}$
when $G_2/K_2$ is a sphere. In fact, the Einstein equation is still valid if we
replace $G_1/K_1$ by any positive Einstein manifold of dimension $d_1=2, 3$ or  $5$.
Note that the dimension of the resulting Ricci-flat manifold is $10$ or $11$.
The {\em explicit} Ricci-flat metrics obtained here belong to a family of complete
Ricci-flat metrics shown to exist earlier in \cite{Bo4} by dynamical systems method.
\end{example}

\begin{example} $($\cite{DW2}$)$ \label{27dim}
One can generalize the first integral in Example \ref{Bergery} for a scalar
curvature function of the same form when $(d_1, d_2)$ satisfy the diophantine
equation $\frac{4}{d_1} + \frac{9}{d_2} =1$ (cf case (4) in Theorem \ref{smallrthm}).
However, the corresponding dimensions can be realized by a principal
orbit only when $(d_1, d_2)=(8, 18)$, so that $n+1=27$. In this case, $G/K$
are certain homogeneous $8$-torus bundles over a product of nine $\C\PP^1$s.
While the Ricci-flat equations can again be solved explicitly, there are no
complete metrics among the solutions. Incompleteness occurs at the finite end
of the cohomogeneity one manifold.
\end{example}

In view of the above three examples, one might wonder whether more conserved
quantities exist if we extend the class of functions (cf (\ref{intform}))
or weaken the equation (cf (\ref{polygenFI})). A partial result addressing
this issue is

\begin{thm} $($\cite{DW9}, Theorem 1.1$)$ \label{quadint}
Suppose that the scalar curvature function of the principal orbit
$G/K$ is of the form $A_1 \exp(-q_1) + A_2 \exp(-q_2)$ or
$A_1 \exp(-q_2) - A_2 \exp(q_1 -2q_2)$ as above. Let $\mathsf F$
be a function on $T^*({\mathcal C})$  of the form
$$ {\mathsf F}(q,p) = {\sf F}_{20}\, p_1^2  + {\sf F}_{02}\, p_2^2
   + {\sf F}_{11}\, p_1p_2 + {\sf F}_1\, p_1 + {\sf F}_2\, p_2 + {\sf F}_0, $$
where ${\sf F}_{ij}$ and ${\sf F}_i$ are $C^1$ functions in $q$ for $0 \leq i,
j \leq 2.$

Suppose further that $\sf F$ satisfies $\{{\sf F}, {\sf H}\} = 0$
on ${\mathscr V}_{\sf H}$. Then $F= \phi {\sf H} + c $ where $\phi$
is a $C^1$ function in $q$ and $c$ is a constant $($ i.e., $\sf F$ is a
trivial first integral $)$ except in the following situations.

Case $($i$):$   With the convention
   $d_1 \leq d_2$, we have either $d_1=1$ or  $4d_1=d_2(d_1 - 1)$.

Case $($ii$):$  We have $d_1=1$ or  $d_2 =2$ or $9d_1=d_2(d_1 -4)$.

Moreover, a non-trivial linear integral, i.e., one with
${\sf F}_{20}={\sf F}_{11}={\sf F}_{02}=0$ and ${\sf F}_1, {\sf F}_2$ not both zero,
exists iff we are in Case $($i$)$ and $d_1=d_2=5.$
\end{thm}

\begin{rmk} \label{rm-quadint}
The first integrals that occur in Theorem \ref{quadint} differ from those
described in Examples \ref{Bergery} - \ref{27dim} essentially by an element in the
ideal generated by the Hamiltonian. Further details about these first integrals as
well as why the case $d_2=2$ in (ii) is special can be found in \cite{DW9}.
\end{rmk}

\smallskip

\noindent{\bf 5C. Painlev{\'e}-Kowalewski analysis}

\smallskip

Besides searching for superpotentials or conserved quantities of the
cohomogeneity one Einstein equations, there is another classical technique for
singling out special cases of these equations which have ``nice" properties
such as integrability. It originated from Kowalewski's prize-winning work on
integrable tops \cite{Kow} and the classification by Painlev{\'e} of second order
ODEs of the form $\ddot{x} = F(\dot{x}, x, t)$ whose solutions have poles as
the only movable singularities (the {\em Painlev{\'e} property}). Recall that a
singularity of a solution of an ODE (with algebraic functions as coefficients)
is called {\em movable} if its location depends on integration constants. The
simplest example is the pole of $x(t)=(t-a)^{-1}$, which is a solution to
$\dot{x} + x^2 = 0$. Movable singularities of solutions are particular to
nonlinear ODEs since singularities of solutions to linear ODEs depend only on
the location of singularities in its coefficients ({\em the fixed singularities}).
The Painlev{\'e} property has been heuristically associated with integrability
since Kowalewski's work.

A modern formulation of Painlev{\'e} analysis for ODEs has been given by
Ablowitz, Ramani, and Segur \cite{ARS}. This involves looking for {\em maximal}
families of solutions of the ODE system in which each dependent variable
is a convergent Laurent series about some common centre (the movable singular point).
One requires at least one of the Laurent series to have a pole, i.e., the
corresponding dependent variable blows up as one approaches the centre.
For an $N$-dimensional system of ODEs, the most general solution should depend
on $N$ arbitrary parameters, so maximality means that the Laurent series
solutions should depend on $N$ parameters with one of these being the location
of the pole. The signal for integrability is strongest if {\em each} Laurent series
blows up and depends on the full number of parameters.

In looking for these Laurent series solutions, one proceeds by first determining
the leading powers (hoping for negative integer powers) and then solving for
the higher order terms by recursion. The recursion operator will fail to be
invertible at various steps which are called the {\em resonances}. At each such
step, there are {\em compatibility conditions} which ensure that the right-hand
side of the recursion relation lies in the range of the recursion operator.
If these are satisfied, then free parameters enter the solution series being
constructed. Finally one needs to prove convergence of the solution series
in a deleted neighbourhood of the singularity.

It may be the case that the number of parameters in the solution series falls
short of the maximal number. Experience has shown that the equation may still have
nice properties if the number of parameters is large compared to $N$.
When one fails to obtain Painlev{\'e} series solutions of the above type, one may
consider Laurent series expansions in some rational power of the independent
variable.

In \cite{DW4}-\cite{DW8} Painlev{\'e} analysis has been applied to the cohomogeneity
one Einstein equations in various cases where Ricci-flat metrics, particularly
ones with special holonomy, have been discovered. The principal orbits $G/K$ are
assumed to have multiplicity free isotropy representation. Using ideas in \cite{AdvM},
we can embed the cohomogeneity one Einstein equations in a quadratic system
(cf \cite{DW6}, \S 2) which may be regarded as describing a Poisson Hamiltonian flow.
The Einstein equations then appear as the restricted flow on a symplectic leaf
cut out by certain constraint equations. In addition to being an interesting fact
on its own, the existence of this embedding has a practical consequence.  Once one
has constructed formal Painlev{\'e} expansion solutions to the quadratic system,
a general majorization argument (cf \cite{DW4}, \S 6) establishes their convergence
in a deleted neighbourhood of the movable singularity once and for all.

Because Painlev{\'e} expansions are meromorphic in some rational power of
the independent variable of the quadratic system, they often represent Ricci-flat
manifolds whose orbit space is $[t_*,\infty)$ with $t_*$ large. The metrics
are complete at infinity, however one generally does not know if they
can be extended beyond $t_*$ and compactified (by adding a lower-dimensional
orbit) to give a global complete smooth metric. Painlev{\'e} analysis seems
particular suited for studying asymptotic behaviour of Einstein metrics.

\begin{rmk} \label{ALC}
In describing the asymptotics of a cohomogeneity one Ricci-flat metric, physicists
have introduced some convenient terminology which we will adopt. The metric is said
to be {\em asymptotically conical} (AC) if for sufficiently large geodesic distances
$t$, the metric is asymptotic to the Ricci-flat cone over some invariant Einstein
metric on the principal orbit. It is {\em asymptotically locally conical} (ALC) if
the principal orbit is a homogeneous circle bundle over some base $N$ and the metric
has asymptotic form
$$ dt^2 + a (d\theta)^2 + t^2 g_N $$
where $a>0$ is a constant, $\theta$ is the coordinate along the circle, and $g_N$ is Einstein
(see also Proposition 2.15 in \cite{DW8}). The fact that the circle has asymptotically
constant radius is significant in $M$-theory (see \cite{GuSp}).
\end{rmk}

Below we describe a sample of the results of the Painlev{\'e} analyses performed
in \cite{DW4}-\cite{DW8}. Note that two of the parameters that are counted represent
the position of the singularity and an overall homothety. These parameters are of course
geometrically trivial, but they will be included in all discussions below.

\smallskip

{\bf (1) Principal Orbit $\SU(3)/T$:} The isotropy representation consists of
$3$ inequivalent irreducible summands, permuted by the Weyl group of $\SU(3)$, which
is the symmetric group on $3$ objects. The cohomogeneity one Ricci-flat system is
therefore $5$-dimensional, in view of the zero energy constraint. Bryant and Salamon
\cite{BrS} as well as Gibbons, Page, Pope \cite{GiPP} have examined this principal
orbit type for metrics with holonomy ${\rm G}_2$. A superpotential was found in \cite{CGLP4},
and the associated first order system (of dimension $3$ as this system automatically
satisfies the zero energy constraint by Proposition \ref{linint}(a)) was identified as
representing ${\rm G}_2$ holonomy metrics (see also Theorem \ref{nonnullthm}(5)).

Painlev{\'e} analysis of the Ricci-flat system yields three families of
meromorphic solutions. There is, up to the action of the Weyl group, a
$5$-parameter family of Painlev{\'e} expansions where only one of the
dependent variables blow up. It corresponds to incomplete Ricci-flat
metrics with generic holonomy. The second family depends on only one
parameter, the position of the singularity. It represents the Ricci-flat
cone metric on one of the $3$ isometric $\SU(3)$-invariant K\"ahler-Einstein
metrics on $\SU(3)/T$.

The third family depends on $3$ parameters and represent metrics with
holonomy in ${\rm G}_2$. It is therefore a maximal family with respect to
the first order subsystem associated to the superpotential. The Painlev{\'e}
expansions are meromorphic in $s^{\frac{1}{5}}$, where $s$ denotes the
independent variable of the quadratic system, and all dependent variables
of the quadratic system actually blow up. Asymptotically, the Ricci-flat
metrics behave like the cone on the normal metric on $SU(3)/T^2$, which
is Einstein but non-K\"ahler. The full $3$-parameter family of solutions of
the first order subsystem was obtained also by Cleyton \cite{Cl} and by
\cite{CGLP4}. In fact, the first order system is equivalent to
the $\SU(2)$ Nahm's equation (well-known to be integrable) by a suitable change
of variables. This was noticed independently in \cite{CGLP4} and in \cite{DW8}.

\smallskip

{\bf (2) Principal Orbit $S^3 \times S^3 \approx ((\SU(2) \times \SU(2) \times
\Delta \U(1)) \ltimes  {\Z}_2)/(\Delta \U(1) \times {\Z}_2)$:} In the above
coset representation of $S^3 \times S^3$, the diagonal $\U(1)$ in $G$ is the
diagonally embedded circle in $S^3 \times S^3$ regarded as right translations
on $\Spin(4)$. The $\Z/2\Z$ in $G$ interchanges the two $\SU(2)$ factors. The
subgroup $K$ is specified by the inclusions
$$ \Delta \un(1) \hookrightarrow \Delta \un(1) \oplus \Delta \un(1)
\hookrightarrow (\un(1) \oplus \un(1)) \oplus \Delta \un(1)
\hookrightarrow \su(2) \oplus \su(2) \oplus \Delta \un(1). $$
The isotropy representation consists of $4$ pairwise inequivalent
irreducible summands, so the Ricci-flat system is $8$-dimensional. There are two
known invariant Einstein metrics on $G/K$, the product metric and the Killing form
metric, see Example \ref{spin4}.

The search for cohomogeneity one metrics with ${\rm G}_2$ holonomy and
$S^3 \times S^3$ as principal orbit type was initiated in \cite{BrS} and
\cite{GiPP}. The specific coset representation of $S^3 \times S^3$ given above
is due to Brandhuber, Gomis, Gubser and Gukov \cite{BGGG}, who found the
superpotential (cf Remark \ref{disconncase}) and showed that the associated
first order system represented the ${\rm G}_2$ holonomy condition. They produced
a complete explicit solution of this system modulo homothety, and gave
an argument base on pertubation theory that this solution lies in a $1$-parameter
family of geometrically distinct solutions.

Painlev{\'e} analysis (cf \cite{DW8}) shows that there are only two families
of solutions where all the dependent variables blow up. The Painlev{\'e} expansions
involve rational powers of the independent variable. The first family depends on
$3$ parameters. The corresponding Ricci-flat metrics are complete at
infinity and asymptotic to a cone over an Einstein metric on the principal
orbit. Within this family is a $2$-parameter family of metrics with holonomy
${\rm G}_2$. These are the ones discoverd in \cite{BrS} and \cite{GiPP}.
The second family depends on $4$ parameters. The corresponding Ricci-flat
metrics are complete at infinity and are asymptotic to a circle bundle
over an Einstein cone (ALC asymptotics). Within this family are two $3$-parameter
families of metrics with ${\rm G}_2$ holonomy. (Recall that one of the parameters
is always the position of the singularity.) These families correspond to the
${\rm G}_2$ metrics described by \cite{BGGG}.

There are further families of Painlev{\'e} expansions which depend on
$6$ or $7$ parameters. In these cases, only $2$ (resp. one) of the
dependent variable actually blow up. The corresponding Ricci-flat metrics
are incomplete at both ends.

Finally, there is a $5$-parameter family of Painlev{\'e} expansions
in which three of the dependent variables blow up. The corresponding
Ricci-flat metrics do not have special holonomy. They are, however, complete
at infinity and are asymptotic to a $2$-torus bundle over a cone on
$S^2 \times S^2$.

\smallskip

{\bf (3) Principal Orbit $N_{p,q}=\SU(3)/\U_{p,q}$:} These are the {\em generic}
Aloff-Wallach spaces, i.e., $\{p,q,-(p+q)=m\}$ is different from $\{1, 1, -2\}$
or $\{1, -1, 0\}$. We also assume that $(p, q)=1$ (so the spaces are simply
connected), and fix one representative among all permutations of $\{p,q,-(p+q)\}$.
The isotropy representation then consists of $4$ irreducible summands and so
the Ricci-flat equations are a $7$-dimensional system. There are exactly two
invariant Einstein metrics on each $N_{p, q}$ by \cite{PP1}.

The search for cohomogeneity one metrics with holonomy $\Spin(7)$ was initiated
in \cite{CGLP4} and \cite{KaY} independently. These authors derived a superpotential
for the Ricci-flat system and showed that the associated first order system
represented the $\Spin(7)$ holonomy condition. An isolated asymptotically
conical solution to the first order system was found, and numerical evidence
for asymptotically locally conical solutions was given.

The Painlev{\'e} analysis in \cite{DW7} produced three families of local Ricci-flat
metrics which are asymptotically locally conical and which contain $\Spin(7)$ holonomy
metrics. The first family depends on $5$ parameters and contains a $4$-parameter
subfamily consisting of metrics with $\Spin(7)$ holonomy. Notice that this is
a full family with respect to the associated first order system. As well, all the
dependent variables in the quadratic system actually blow up. The corresponding
metrics are asymptotic to a circle bundle over a cone on $\SU(3)/T$ equipped with
the Killing form metric. The circle fibres have asymptotically constant radii.
The numerical examples indicated in \cite{CGLP4} and \cite{KaY} belong to this family.

There is a second family of Painlev{\'e} expansions which depends on $6$ parameters
and contains a $3$-parameter family of $\Spin(7)$ metrics. In this family, only two
of the dependent variables blow up. The third family depends on
$7$ parameters and contain a $4$-parameter family of $\Spin(7)$ metrics. So
we again have a Painlev{\'e} family depending on the full number of parameters.
In this case, however, only one of the dependent variables blow up.

Painlev{\'e} analysis also shows that there are families of local Ricci-flat
metrics which are asymptotically conical, i.e., asymptotic to the cone over one
of the two invariant Einstein metrics of $N_{p, q}$. However, except for finitely many
pairs $(p, q)$, these families depend on at most $2$ parameters. This is because
the condition of having non-trivial {\em rational resonances} implies that there is a
rational point on a certain smooth curve $y^2 = p(x)$ where $p(x)$ is an explicit
polynomial of degree $6$. Finiteness then follows from Faltings' solution of the
Mordell conjecture and the relation between $x$ and the explicit expression \cite{CaR}
of the invariant Einstein metrics on $N_{p, q}$ in terms of $p$ and $q$. So the
Painlev{\'e} analysis is consistent with the sparseness of the known asymptotically
conical solutions.

\smallskip

{\bf (4) Principal Orbit $S^{4m+3} \approx \Sp(m+1)\U(1)/\Sp(m)\Delta\U(1)$:}
The isotropy representation consists of $3$ irreducible summands of dimensions
$1, 2$ and $4$. So the Ricci-flat equations are a $5$-dimensional system. There
are two invariant Einstein metrics on the principal orbit--the constant curvature
metric and the Jensen metric \cite{Jen2}.

When $m=1$, the principal orbit is $7$-dimensional. The search for cohomogeneity one
metrics with $\Spin(7)$ holonomy and this principal orbit type was first undertaken
in \cite{BrS} and \cite{GiPP}. These authors found a complete $\Spin(7)$ holonomy
metric which is asymptotic to the cone over the Jensen metric. A superpotential
for the Ricci-flat system was found in \cite{CGLP3} (cf Theorem \ref{nonnullthm}(4)),
and the associated first order system was shown to represent the $\Spin(7)$ holonomy
condition. This system turns out to be integrable, and one obtains a $3$-parameter
family of solutions (of which $2$-parameters are geometrically trivial).

For general $m$, Painlev{\'e} analysis yields two $3$-parameter families of
solutions of the quadratic system in which all the dependent variables blow up.
The corresponding Ricci-flat metrics have generic holonomy and are complete
at infinity. In one family, the metrics are asymptotic to the cone over the constant
curvature metric. In the other family, the metrics are asymptotic to a circle bundle
over a cone on $\C\PP^{2m+1}$ where the circle fibres have asymptotically constant
radii.

When $m=1$, there are two further families of Painlev{\'e} expansions in which
all dependent variables blow up. One of these is a $4$-parameter family meromorphic
in $s^{1/4}$, where $s$ is again the independent variable of the quadratic
system. The Ricci-flat metrics are complete at infinity and are asymptotic to
the cone on the Jensen metric. Within this family is a $2$-parameter family
of Painlev{\'e} expansions representing metrics with $\Spin(7)$ holonomy. These
are the ones found by Bryant-Salamon and Gibbons-Page-Pope. The other family
also depends on $4$ parameters. It is meromorphic in $s^{1/5}$. The
corresponding Ricci-flat metrics are complete at infinity and are asymptotic to
a circle bundle over a cone on $\C\PP^{2m+1}$. Within this family is a $3$-parameter
family of Painlev{\'e} expansions representing metrics with $\Spin(7)$ holonomy.
These are exactly the metrics found in \cite{CGLP3}.

Are there analogues of the last two Painlev{\'e} families when $m > 1$ ? It turns out
that during the Painlev{\'e} analysis one reaches a stage at which the resonance should
be a rational root of an equation of the form $x(x-1)=R(m)$ where $R(m)$ is a rational
expression in $m$. When $m=1$ this equation has a rational solution and so
the recursion can proceed, giving rise to the above families. When $m >1$, after some
further analysis, one sees that the above equation implies that $m$ is an integral
point on a certain elliptic curve. Besides $m=1$, no further integral points were found
using the computer programs RATPOINTS and SIMATH in the range $1 < m < 10^6$. So it
appears that Painlev{\'e} analysis can detect the presence of special holonomy solutions.

\section{\bf Complete Einstein Metrics on Fibre Bundles}

In this section we will discuss the construction of Einstein metrics
on fibre bundles using modifications of the Kaluza-Klein framework. This set-up is
well-known in physics, and the mathematical formulation can, for example, be found
in \cite{Be} (particularly 9.61), \cite{Wa2}, or \cite{Wa3}. In order to fix notation
and clarify some issues, we will briefly review the set-up below.

Let $\mathcal P$ be a smooth principal $H$-bundle over a manifold $N$ of dimension
$m$ where $H$ is a compact Lie group. Suppose that $H$ acts smoothly and almost
effectively on another manifold $F$ of dimension $d$. Let $M$ denote the total space
of the associated fibre bundle ${\mathcal P} \times_H F$ and $\pi: M \rightarrow N$
denote the projection map. The basic input data consists of a metric $g^*$ on $N$,
a principal $G$-connection $\omega$ on $\mathcal P$ with curvature form $\Omega$, and
an $H$-invariant metric $h$ on $F$. From this data, one can construct a unique metric
${g}$ on $M$ such that $\pi$ is a Riemannian submersion with totally geodesic
fibres.

In order to write down the Einstein condition for ${g}$, we need to introduce
two inclusion maps. For a point $[p,x] \in {M}$, we have the inclusion
of the fibre
$$ i_p: F \hookrightarrow {M}, \,\,\,\, i_p(x) = [p, x]$$
which satisfies $i_{pa} = i_p \circ a$ for $a \in H$. We also have the inclusion
$$ j_x: {\mathcal P} \hookrightarrow {M}, \,\,\,\,\, j_x(p) = [p, x]$$
which satisfies $j_{ax} = j_x \circ R_a$, where $R_a$ denotes the right action
of $H$ on $\mathcal P$. If $V \in \g$, then it induces a vector field $\overline{V}$
on $F$ which is Killing for the metric $h$. Except when $F=H$ and $H$ acts on itself
by left translation (i.e., ${M} \approx {\mathcal P}$), $\overline{V}$ may
have zeros and so cannot be used as part of a global frame on $F$. We also note
the basic relationship $a_*(\overline{X}_{x}) =(\overline{\Ad_a X})_{ax}$.

\begin{prop} \label{Kaluza} The Einstein condition for the metric ${g}$ is
given by

$($i$)$ the connection $\omega$ is Yang-Mills,

$($ii$)$ for all vertical tangent vectors $U, V$ at $[p,x] \in {M}$ we have
\begin{equation} \label{fibreq}
\Ric(h)(i_{p*}^{-1}(U), i_{p*}^{-1}(V)) + \frac{1}{4} \sum_{i,j} \,
    h( \overline{{\Omega({\tilde e}_i, {\tilde e}_j)}}_x, i_{p*}^{-1}(U))
    \  h(\overline{{\Omega({\tilde e}_i, {\tilde e_j})}}_x, i_{p*}^{-1}(V))
  = \Lambda \, h(i_{p*}^{-1}(U), i_{p*}^{-1}(V)),
\end{equation}

$($iii$)$ for all horizontal vectors $X, Y$ at $[p,x] \in {M}$ we have
\begin{equation} \label{baseq}
\Ric(g^*)({\pi_*}(X), {\pi_*}(Y)) -\frac{1}{2} \,\sum_i  \, h(\overline{{\Omega( \tilde{X},
{\tilde e}_i)}}_x, \overline{{\Omega((\tilde{Y}, {\tilde e}_i)}}_x)
= \Lambda \, g^*({\pi_*}(X), {\pi_*}(Y)),
\end{equation}
where $\Lambda$ is the Einstein constant, $\tilde{X}, \tilde{Y}$ are respectively
the $\omega$-horizontal lifts of $\pi_*(X), \pi_*(Y)$ to
$p \in {\mathcal P}$, $\{e_1, \cdots, e_m\}$ is an
arbitrary orthonormal frame at $\pi([p, x])$ in $N$, and $\{\tilde{e}_1, \cdots, \tilde{e}_m \}$
is its horizontal lift to $p \in {\mathcal P}$.
\end{prop}

Note that unless ${M}$ itself is the principal bundle $\mathcal P$ {\em and}
$h$ is a bi-invariant metric on $H$, Equation (\ref{baseq}) is not an equation on the
base. The Einstein condition above is actually very restrictive. By 9.62 in \cite{Be},
$g^*$ and $h$ must have constant scalar curvature and the pointwise norm $|{\Omega}|$
must be everywhere constant. Moreover, the second term on the left of Equation (\ref{baseq})
must be the same everywhere on a fixed fibre, and the second term on the left of
Equation (\ref{fibreq}) should be independent of $p$.

One situation in which the above system can be analysed is that of fibrations
of homogeneous spaces. In this case all the difficulties mentioned above disappear,
leaving behind a more tractable set of conditions. Explicitly, in the fibration
$$ H/K \longrightarrow G/K \longrightarrow G/H$$
where $K, H$ are closed subgroups of a compact Lie group $G$, we have
$G/K \approx G \times_H (H/K)$. So $\mathcal P$ is just $G$, with $F = H/K$.
A choice of Yang-Mills connection is nothing but an $\Ad_H$ invariant
$b$-orthogonal decomposition $\g = \h \oplus \p_{-}$. Its curvature form
clearly has constant norm. The use of such fibrations to construct invariant
Einstein metrics on the total space $G/K$ was first systematically studied by
G. Jensen \cite{Jen2}. It led to the discovery of the first variable curvature
Einstein metrics on spheres $S^{4m+3}$ as well as numerous left-invariant metrics
on compact Lie groups (see \cite{Jen2} and \cite{DZ}).

Recently, Ara{\'u}jo \cite{Ara1} worked out some consequences of the Einstein
condition for a $G$-invariant metric $g$ on the total space $G/K$ under the assumption
that the isotropy representations of $H/K$ and $G/H$ have no multiplicities
and no common $\Ad_K$ irreducible summands. These necessary conditions are algebraic
in nature and involve Casimir type operators associated to the irreducible
summands in the isotropy representations of $H/K$ and $G/H$. Let $G$ be compact
and semisimple. Then $g$ is said to be {\em binormal} if the fibre metric $h$
and base metric $g^*$ are induced by constant multiples of the Killing form metric
of $G$. In this case, Ara{\'u}jo gives a necessary and sufficient condition for $g$
to be Einstein. The special situation where $h$ and $g^*$ are both Einstein
is further analysed in \cite{Ara1}. There, a non-normal homogeneous Einstein metric is found
on $(G \times \cdots \times G)/\Delta G$ when the number of factors of $G$
is $\geq 5$. (The Killing form metric is known to be Einstein.) In \cite{Ara2}
a partial classification is given under the further assumption that $H/K$
is a simply connected symmetric space and $G/H$ is an irreducible symmetric space.

\smallskip

To inject some flexibility into Equations (\ref{fibreq}) and (\ref{baseq}), a
simple modification is to assume that $H$ acts on $F$ with cohomogeneity one and
simultaneously multiply the base metric $g^*$ by an appropriate function of the
orbit space coordinate $t$. In other words, we consider on ${M}$ a metric of the form
$$ g = dt^2 + h_t + \alpha(t)^2 \pi^* g^*$$
(in which the connection $\omega$ is suppressed). More generally, one can
replace $g^*$ by a one-parameter family of metrics on $N$ depending smoothly
on $t$. Notice that when $N$ is itself homogeneous, e.g., $N=G/(A\cdot H)$,
then we are back to the cohomogeneity one set-up with ${M}=
(G/A)\times_H F$. If the principal orbit type in $F$ is $H/L$, then the
principal orbit type in ${M}$ is
$(G/A) \times_H  (H/L) \approx G/(A\cdot L).$ On the other hand, $N$ need
not have any symmetries, in which case $M$ will in general also
have little symmetry.

The basic geometry of the above metric ansatz is very similar to that of the
cohomogeneity one case. The role of the principal orbit is played by
the hypersurface ${\mathcal Z}:={\mathcal P} \times_H (H/L)$. After removing
from the manifold ${M}$ the submanifolds of the form
${\mathcal P} \times_H (H/K)$ where $H/K$ is a singular orbit in $F$, we again
have an equidistant family of hypersurfaces $I \times {\mathcal Z}$. As seen
in \cite{EW}, the Einstein equations are again of the form (\ref{eq1})-(\ref{eq4}).
For each fixed $t$, the metric $h_t + \pi^* g^*_t$ on $\mathcal Z$ makes $\pi$
into a Riemannian submersion with totally geodesic fibres. Many features of
the cohomogeneity one set-up carry over if we further choose $g^*_t$ and $\omega$
such that the hypersurfaces $\{t\} \times {\mathcal Z}$ have constant mean
curvature and constant scalar curvature. As in the cohomogeneity one case,
there are additional conditions at the boundaries of $I$, the orbit space of
$F$, to ensure that the metric $g$ is smooth and complete.

\smallskip

\noindent{\bf 6A. Bundles with abelian structural group}

\smallskip

The above modification of the Kaluza-Klein set-up has been used to construct
many Einstein metrics on fibre bundles, as was already reported in \S 3 of
\cite{Wa3}. A very fruitful direction has been to take $\mathcal P$ to be a
circle bundle over a product of Fano K\"ahler-Einstein manifolds such that the
Euler class is a rational linear combination of the anti-canonical classes of
the base factors. The fibre $F$ is $\C$, $S^2$ or $\R\PP^2$ with the standard
action of the circle by rotation. The family $g^*_t$ consists of products of
different scalings of the K\"ahler-Einstein metrics, and the Yang-Mills connection
$\omega$ is just the connection on the circle bundle whose curvature form is
$-2\pi i$ times the harmonic representative of its Euler class. Notice that the
connection remains Yang-Mills as the base metrics vary with $t$. For fixed $t$,
the curvature form has constant norm, and the hypersurfaces indeed have constant
scalar curvature and mean curvature. The resulting manifolds $M$ have very
little symmetry if the base has little symmetry.

The first people to exploit the above set-up were Calabi, Page, and
B{\'e}rard Bergery. They treated the case where the base consists of a single
Fano KE factor. Calabi was interested in non-positive K\"ahler Einstein metrics
on holomorphic line bundles. Page \cite{P} constructed the first compact example
of a non-homogeneous Einstein metric, and B{\'e}rard Bergery \cite{BB} extended
Page's construction to higher dimensional and non-compact cases after first
reformulating his work in the cohomogeneity one setting (cf also \cite{PP2}).
A few years later, Koiso and Sakane \cite{KoS1}, \cite{KoS2} used the same set-up
to construct the first examples of inhomogeneous positive K\"ahler-Einstein
metrics on certain projectivized holomorphic bundles over products of Fano KE
manifolds or coadjoint orbits.

Notice that if we let the base $N$ to be a coadjoint orbit, then $g^*_t$ can be
allowed to vary among the homogeneous K\"ahler metrics on $N$ (even for different
invariant complex structures) without having to be Einstein. This allows for more
examples, but then we are back to the cohomogeneity one set-up. In \cite{DW1} a
classification is given of non-positive K\"ahler-Einstein metrics of cohomogeneity
one under the assumptions that $G$ is compact connected and semisimple, and the
principal orbit has multiplicity free isotropy representation. The latter is a
generic condition, as ``most" circle bundles over a coadjoint orbit have multiplicity
free isotropy representation. (The situation of the Aloff-Wallach spaces is a good
example.) The classification shows that if $M$ is in addition neither
reducible nor hyperk\"ahler, then it is a line bundle over a coadjoint
orbit or a blow-down of it (under suitable conditions) along its zero section.

For the positive case, the classification is addressed in both \cite{DW1} and
\cite{PoS1}. The assumptions in \cite{DW1} are the same as in the non-positive
case above, while the assumptions in \cite{PoS1} are somewhat weaker. The upshot
in both papers is that the singular orbits are also coadjoint orbits and all
K\"ahler-Einstein metrics are obtained via the Koiso-Sakane construction.
If one turns away from the assumptions in these papers, then it is possible
for the cohomogeneity one K\"ahler manifold to have exactly one complex
singular orbit. In this situation, there are further examples of K\"ahler-Einstein
metrics. The first constructions were given in \cite{GuCh}. In \cite{PoS2}, a
classification of cohomogeneity one positive K\"ahler-Einstein manifolds was given
under the assumptions that (i) there is a unique complex singular orbit of complex
codimension $1$ and (ii) the principal orbits are Levi-nondegenerate. Many new
K\"ahler-Einstein metrics result from this classification, including ones on three
infinite families of K\"ahler manifolds.

We turn next to the non-K\"ahler case.

Examples of hermitian but non-K\"ahler Einstein metrics on certain $S^2$ and
$\R\PP^2$ bundles over a product of Fano KE manifolds were constructed in the
thesis of J. Wang (cf \cite{WaW}). The $\R\PP^2$ case is a generalization
of the work of Page and B{\'e}rard Bergery, and by lifting to the corresponding
$S^2$ bundle one gets a different Einstein metric than the one constructed directly
on it.  We note that it often happens that these Einstein metrics exist when the
Futaki invariants of the $S^2$ bundles are nonzero. Non-positive non-K\"ahler Einstein
metrics were also constructed on certain complex line bundles over Fano KE products.
A blow-down analysis similar to that in \cite{DW1} can be applied to the above bundle
type metrics to obtain further examples. These subsequent results are due to Dancer and
myself, and are stated as Theorems 3.3-3.6 in \cite{Wa3}. It should be noted that the
Ricci-flat metrics constructed in \cite{WaW} and \cite{Wa3} were rediscovered in
\cite{CGLP6}. In fact the Ricci-flat equations admit a superpotential \cite{CGLP6}
(see also Theorem \ref{nullthm}) which singles out the Calabi-Yau metrics. Painlev{\'e}
analysis of this Ricci-flat system can be found in \cite{DW7}.

\smallskip

Circle bundles over Fano KE products are also useful for constructing explicit examples
of {\em conformally compact Einstein} manifolds as they bound the corresponding disc
bundles. Let $\overline{M}$ be a smooth manifold with boundary and set $M=\overline{M}
\setminus \partial M$. Recall that a complete metric $g$ on $M$ is said to be conformally
compact Einstein if $g$ is a negative Einstein metric and there exists a defining function
$\rho:\overline{M} \rightarrow \R_{+} \cup \{0\}$ such that $\overline{g} = \rho^2 g$
extends to a $C^2$ metric on $\overline{M}$. The {\em conformal infinity} of $(M, g)$ is
then $\partial M$ together with the conformal class of $\overline{g}\,|\partial M$.
For more information regarding conformally compact Einstein manifolds, see for
example \cite{An}.

Now let ${\mathcal P}_q$ be a principal $S^1$ bundle over a product
$(N_1, h_1) \times \cdots \times (N_r, h_r)$ of Fano KE manifolds with Euler class
$ q=q_1\alpha_1 + \cdots + q_r \alpha_r$ where $q_j$ are nonzero integers,
$c_1(N_j) = p_j \alpha_j$ with $p_j \in \Z_{+}$ and
$\alpha_j \in H^2(N_j; \Z)$ indivisible. Assume that the complex dimension
of $N_j$ is $n_j$ and the KE metric $h_j$ is normalized so that the first Chern
number $p_j$ is the Einstein constant. In addition let $(X, \gamma)$ denote an
arbitrary closed Einstein manifold with Einstein constant $\mu$ and dimension $k$.
(The possibility that $X$ reduces to a point is included.) Then we have a manifold
$({\mathcal P}_q \times_{S^1} {\mathcal D}) \times X$ with boundary
${\mathcal P}_q \times X$, where $\mathcal D$ is a $2$-disc of some finite radius.
Let $\overline{g}$  be a bundle-type metric on $\overline{M}:={\mathcal P}_q \times_{S^1}
{\mathcal D}$ of the form
$$ \overline{g}:=dt^2 + f(t)^2 d{\theta}^2 + g_1(t)^2 h_1 + \cdots + g_r(t)^2 h_r$$
where $f(t), g_j(t)$ are non-negative functions and the connection $\omega$ on
the circle bundle used to define the metric has been suppressed. (The
curvature form of $\omega$ is, as usual, $-2\pi i$ times the harmonic representative
of the Euler class $q$.) One now seeks a defining function $\rho(t):\overline{M} \rightarrow
\R_{+}$ such that $\rho^{-2}(\overline{g} \oplus \gamma)$ is a complete negative
Einstein metric on $M \times X$ or a suitable blow-down of it along
the zero section of $\overline{M}$. In \cite{Chn1} it is shown that it suffices to
solve the {\em quasi-Einstein type equation}
\begin{equation} \label{quasiE}
\Ric(\overline{g}) + (d-2)\, \frac{{\rm Hess}_{\overline{g}}(\rho)}{\rho} = \mu \overline{g}
\end{equation}
on $\overline{M}$, where $d=k+\dim \overline{M}= k+ 2(1+ \sum n_j)$.

One can then use methods similar to those in \cite{WaW}, \cite{DW1}, and \cite{DW14}
to obtain explicit solutions of this equation.

\begin{thm} $($\cite{Chn1}, \cite{MaS}$)$ \label{confcomp} With notation as above,
assume further that $r \geq 2$, $(N_1, h_1)$ is $\C\PP^{n_1}$ with the Fubini-Study
metric $($suitably normalized$)$, and $|q_1|=1$.

$($1$)$ Suppose $\mu > 0$. Then there is a $1$-parameter family of non-homothetic
  conformally compact Einstein metrics on $E \times X$ where $E$ is the rank $n_1 +1$
  complex disc bundle over $N_2 \times \cdots \times N_r$ resulting from blowing down
  the zero section of $\overline{M}$ along $N_1$.

$($2$)$ Suppose $\mu \leq 0$. If $(n_1+1)|q_j| > p_j$ for all $2 \leq j \leq r$, then
  the same conclusion as in $($1$)$ holds.
\end{thm}

\begin{rmk} \label{chen1}
$($a$)$ When $n_1=0$ in the above, it means there is no blow down, i.e., $E=M$.
    When $X$ reduces to a point, we can still solve Equation (\ref{quasiE}) with $\mu$
    having the indicated sign.

$($b$)$ In the case where $\mu < 0$, Theorem \ref{confcomp} still holds if
     the base factors $N_2, \cdots, N_r$ are just K\"ahler-Einstein. If $N_j$
     is Ricci-flat, i.e., $p_j=0$, we must in addition assume that the K\"ahler
     class of $h_j$ is $2\pi \alpha_j$. For non-positive KE factors, the
     condition $(n_1+1)|q_j| > p_j$ is of course trivially satisfied.

$($c$)$ The negative Einstein metrics given in Theorem \ref{confcomp}(1) where $X$
      reduces to a point are known (cf \cite{BB}, \cite{C}, \cite{Pe}, \cite{WaW}
      \cite{Wa3}). However, conformal compactness was not explicitly noted
      in these references.

(d) In the special case $n_1 =1, n_2=0, k=0$ we get a $1$-parameter family of
    conformally compact metrics on the $4$-ball with conformal infinity a squashed
    $3$-sphere. These are the AdS-Taub-NUT metrics of Hawking, Hunter and Page \cite{HaHP}.
    Similarly, the case $n_1=0, n_2 =1, k=0$ gives the AdS-Taub-Bolt metrics in \cite{HaHP}.

(e) In addition to the construction of conformally compact Einstein metrics, the paper
   \cite{Chn1} also examines the conformal invariants (renormalized volumes, conformal
   anomalies) associated to the Einstein manifolds. Furthermore, it is shown that whenever
   the total dimension $k+2(1 + \sum n_j)$ is odd, the conformal infinity of each Einstein
   manifold given by Theorem \ref{confcomp} always contains a representative with zero
   $Q$-curvature. The blow-down analysis is also treated carefully.
\end{rmk}

\smallskip

Instead of circle bundles over a Fano KE manifold we can also use $2$-torus bundles to
construct Einstein metrics on associated fibre bundles. We need to fix an isomorphism
$T^2 \approx S^1 \times S^1$ first. Then since $T^2$ acts with cohomogeneity one on
$S^3$ and on the solid torus $D^2 \times S^1 \subset S^3$, these spaces can be taken as
the fibre $F$ in the modified Kaluza-Klein ansatz above. We now take ${\mathcal P}_q$
to be the principal $2$-torus bundle over a Fano KE manifold $N$ with Euler class
$0 \neq q=(q_1 \alpha, q_2 \alpha) \in H^2(N;\Z) \oplus H^2(N; \Z)$, where as before
we write the first Chern class of $N$ as $p\alpha$ with $p>0$ and $\alpha$ indivisible,
and we normalize the K\"ahler-Einstein metric $g^*$ on $N$ to have Einstein constant $p$.
Note that via an automorphism of the torus, ${\mathcal P}_q$ is diffeomorphic to
${\mathcal P}_{q_0} \times S^1$ where $q_0$ is the greatest common divisor of $q_1$ and $q_2$.

\begin{thm} $($\cite{Chn2}$)$ \label{2torusthm}
$($1$)$ Let $M$ denote ${\mathcal P}_q \times_{T^2} (D^2 \times S^1)$ with $|q_2| \neq 0$.
Then there is a $2$-parameter family of conformally compact Einstein metrics on $M$
with conformal infinity ${\mathcal P}_q$.

$($2$)$ There is also a $2$-parameter family of complete Ricci-flat metrics on $M$.
These have sub-Euclidean volume growth and quadratic curvature decay.

$($3$)$ For $|q_1|>|q_2| >0$ there is a positive Einstein metric on the associated
$S^3$ bundle of ${\mathcal P}_q$.
\end{thm}

\begin{rmk} \label{chen2}
(a) If $q_2=0$ in case (1) above, then $M \approx S^1 \times ({\mathcal P}_{q_1} \times_{S^1} D^2)$.
In this case, by Theorem \ref{confcomp}(1) there is a $1$-parameter family of conformally
compact metrics provided that $|q_1| > p.$ (The space $X$ is $S^1$, $N_2=N$ and $n_1=0$.)

(b) The special case of Theorem \ref{2torusthm}(3) where $N=S^2$ was obtained earlier
 in \cite{HaSY} (cf Theorem 1 there). Because of the low dimension, we get only two
 diffeomorphism types of the total space, and hence one gets infinitely many Einstein
 metrics on the two $S^3$ bundles over $S^2$. Some geometric properties of these
 Einstein manifolds (volume, spectrum of the Laplace-Beltrami operator and of the
 Lichnerowicz Laplacian, geodesics) have been studied in \cite{GiHY}.

(c) The topological and smooth classifications of the $S^3$ bundles over $\C\PP^2$
 in Theorem \ref{2torusthm}(3) has been studied in \cite{Chn2} (see Theorem 1.13 there)
 using the invariants of Kreck and Stolz. There are in particular infinitely many
 pairs $(q_1, q_2)$ and $(\tilde{q}_1, \tilde{q}_2)$ such that the corresponding
 $3$-sphere bundles are diffeomorphic, as well as infinitely pairs such that the
 bundles are homeomorphic but not diffeomorphic.

(d) The Einstein metrics in Theorem \ref{2torusthm}(3) have finite isometries coming
from the right action on $S^3$ by diagonally embedded cyclic subgroups of $T^2$.
Hence one also obtain Einstein metrics on the corresponding lens space bundles by
taking quotients.

(e) Each conformal infinity occurring in the examples from Theorem \ref{2torusthm}(1)
admit a representative with zero $Q$-curvature.
\end{rmk}

The paper \cite{HaSY} contains another construction of Einstein metrics on sphere
bundles over $S^2$. This was generalised in \cite{LuPP} as follows.

\begin{thm} $($\cite{LuPP}$)$ \label{spherebundle}
Denote by ${\mathcal P}_q$ the principal circle bundle with Euler class $q\alpha$ over
a Fano KE manifold $N$, where $c_1(N) = p \alpha$ with $p>0$ and $\alpha$ indivisible.
Let $\mathcal Q$ be the bundle ${\mathcal P}_q \times \SO(m+1)$ where $m>1$. Consider
the $S^{m+2}$ bundle over $N$ associated to $\mathcal Q$ via the cohomogeneity one linear
action of $\U(1) \times \SO(m+1) \approx \SO(2) \times \SO(m+1)$ on $F=S^{m+2} \subset \R^{m+3}$.
If $0<q<p$, then there is a positive Einstein metric on the total space of this sphere
bundle.
\end{thm}

Note that the principal orbits in $F$ are $\SO(2) \times (\SO(m+1)/\SO(m))$
and the singular orbits are $\SO(2)$ and $S^m$. The Einstein equations in the
above situation are very similar to those in \cite{BB} and \cite{PP2}, and
can be integrated explicitly by the same method.

In the above construction, if we replace $F$ by $S^2 \times S^m$, i.e., we collapse
the $\SO(2)$ factor at both endpoints of the orbit space, then one obtains
a warped product positive Einstein metric on $({\mathcal P}_q \times_{S^1} S^2) \times S^m$
when $0 < q < p$. Finally, complete Ricci-flat metrics can also be constructed in
two cases. In the first case, one lets $F = \R^2 \times S^m$, i.e., we delete the
point at infinity from $S^2 = \C\PP^1$ in the previous compact case. The result, assuming
again that $0<q<p$, is a complete warped product Ricci-flat metric over the complex line
bundle ${\mathcal P}_q \times_{S^1} \C$. In the second case, one takes $N=\C\PP^k, \, k\geq 1$,
$q=1$, and $F = \R^2 \times S^m$. We can now blow down the zero section of the line
bundle to get $\R^{2k+2} \times S^m$. A complete Ricci-flat metric of warped product
type can be constructed on this space. Further details can be found in \cite{LuPP}.

A more complicated generalisation of the isolated example in \cite{HaSY} consists of
using a $k$-torus bundle ${\mathcal P}_q$ over a Fano KE manifold $N$ with $k \geq 2$
\cite{GLPP}. One fixes an isomorphism of $T^k$ with a $k$-fold product of circles.
Then the Euler class of such a bundle is given by $(q_1\alpha, \cdots, q_k \alpha)
\in H^2(N; \Z)^{\oplus k}$ where $q_j \in \Z$ and $c_1(N) = p \alpha$ as usual. For each
choice of $k$ integers $(a_1, \cdots, a_k)$, there is an action of $T^k$ on $\R^{2k}$
or $\R^{2k+1}$ where $T^k$ is viewed as the usual maximal torus in $\SO(2k)$ (resp. $\SO(2k+1)$)
and the $j^{\rm th}$ circle factor acts through rotation by $a_j \theta_j$. The fibre
$F$ is then $S^{2k-1}$ (resp. $S^{2k}$), which is of cohomogeneity $k-1$ (resp. $k$)
for the torus action. In \cite{GLPP} reasonably strong numerical evidence is given of
a positive Einstein metric on ${\mathcal P}_q \times_{T^k} S^{2k-1}$ whenever all
$q_j > 0$ and $N=S^2$. There are explicit exact solutions when all $q_j$ are equal.
Notice that when $k=2$ we get back the Einstein metrics in \cite{HaSY}, while
Theorem \ref{2torusthm}(3) gives explicit solutions precisely in the case when $k=2$
and the base is replaced by an arbitrary Fano KE manifold. Smooth Einstein metrics exist on
${\mathcal P}_q \times_{T^k} S^{2k}$ only when (up to permutation) $q_2 = \cdots = q_k =0$.
This means that we are back in the situation of two paragraphs ago.

\smallskip

\noindent{\bf 6B. Complete cohomogeneity one metrics with ${\rm G}_2$ holonomy}

\smallskip

Another class of examples of bundle-type Einstein metrics where the structural group
of the bundle is non-abelian consists of complete metrics of cohomogeneity one with
holonomy ${\rm G}_2$. The first examples are due to Bryant-Salamon \cite{BrS} and
Gibbons-Page-Pope \cite{GiPP}. These metrics are of modified Kaluza-Klein type
and occur on the bundle of anti-self-dual $2$-forms on $S^4$ and $\C\PP^2$.
In the first case, as a cohomogeneity one manifold, the principal orbit is
$\C\PP^3 = \Sp(2)/\Sp(1)\U(1)$ and the singular orbit is $\Sp(2)/\Sp(1)\Sp(1) \approx S^4$.
In the second case, the principal orbit is the coadjoint orbit $\SU(3)/T^2$ and
the singular orbit is $\SU(3)/({\rm S}\U(2)\U(1)) = \C\PP^2$. From the bundle perspective,
as $S^4$ and $\C\PP^2$ are the only two self-dual positive Einstein $4$-manifolds
(\cite{Hi1}), the principal orbits are the corresponding twistor spaces.
The connection that is used to construct the ${\rm G}_2$ metrics is the well-known
Yang-Mills connection of the associated principal $\SO(3)$  bundle of the twistor bundle.
A third example due to the same two teams is a ${\rm G}_2$ metric on the real
spinor bundle of $S^3$. Since the $3$-sphere is parallelizable, this bundle is
topologically $S^3 \times \R^4$.

The possible principal orbit types of a cohomogeneity one ${\rm G}_2$ manifold were
determined in \cite{ClS}. There are three cases when the group $G$ is simple:
${\rm G}_2/\SU(3)$, $\Sp(2)/\Sp(1)\U(1)$, and $\SU(3)/T^2$. For these principal
orbits, Cleyton and Swann classified all the complete ${\rm G}_2$ holonomy metrics.
These turn out to be respectively the flat metric on $\R^7$, and the first two
examples in the above paragraph. They also showed that the ${\rm G}_2$ equations
for the principal orbit type $\SU(3)/T^2$ have a non-trivial $1$-parameter family
of solutions, only one member of which is smooth and complete.  The $1$-parameter
family was also discovered independently but slightly later in \cite{CGLP4} and
\cite{DW8}.

The principal orbit types for which the group $G$ is non-simple are $S^3 \times S^3$,
$(\SU(3)\U(1))/\SU(2)$, $\SU(2)\times T^3$ and $T^6$. The case of $G=S^3 \times S^3 = \Spin(4) $
has been studied in \cite{BGGG} and \cite{CGLP4}. The homogeneous geometry of this
group manifold is rather complicated (see Example \ref{spin4}). In \cite{BGGG}, by
enlarging $G$ in a clever way, one arrives at a coset description for $S^3 \times S^3$
in which the isotropy representation consists of only four inequivalent irreducible summands.
The ${\rm G}_2$ ODEs become more manageable, and as a result, an explicit complete
solution with ALC asymptotics was constructed. In fact, infinity looks like a circle
bundle over a cone over $S^2 \times S^3$. Numerical evidence was also given
for a $1$-parameter family of complete solutions. The original example of \cite{BrS}
and \cite{GiPP} is a solution of a $2$-dimensional subsystem of the $4$-dimensional
system and has AC asymptotics.

\begin{rmk} \label{apostolov}
A construction of ${\rm G}_2$-holonomy metrics via $2$-torus bundles over a
hyperk\"ahler $4$-manifold was given in \cite{GLPS} and generalized in \cite{ApS}.
This construction is further exploited in \cite{GS}. We refer the reader to
\cite{ApS} for the more general local description of the construction. Here we give
a brief account of the global case. Suppose $(N, g^*)$ is a hyperk\"ahler $4$-manifold
with associated K\"ahler forms $\omega_1, \omega_2$ and $\omega_3$. Suppose further
that there is a form $\omega_0$ of type $(1, 1)$ which defines an almost complex
structure on $N$ inducing the opposite orientation to that of the hyperk\"ahler structure.
Let $q$ and $s$ be two real parameters, and assume that the de Rham cohomology classes
$\frac{1}{2\pi}[q \omega_0 + s \omega_1]$ and $\frac{1}{2\pi}[\omega_2]$ are integral.
Let $\mathcal P$ be the principal $2$-torus bundle classified by the above two cohomology
classes. Then there is a ${\rm G}_2$-holonomy metric on $(a, b) \times {\mathcal P}$ for
some suitable open interval $(a, b)$.
\end{rmk}

\smallskip

\noindent{\bf 6C. Complete cohomogeneity one metrics with ${\rm Spin}(7)$ holonomy}

\smallskip

We will finally describe efforts to construct complete metrics with $\Spin(7)$
holonomy on vector bundles. The metrics will be of cohomogeneity one type, but one
can also view them as examples of the modified Kaluza-Klein ansatz where the
structural group is non-abelian.

 The earliest complete example is due to Bryant-Salamon \cite{BrS} and
Gibbons-Page-Pope \cite{GiPP}. The underlying manifold is the negative spinor
bundle over $S^4$. As a cohomogeneity one manifold, the principal orbit
is $G/K=(\Sp(2)\times \Sp(1))/(\Sp(1)\Delta \Sp(1))$ and the singular orbit
is $G/H=(\Sp(2)\times \Sp(1))/(\Sp(1)\Sp(1)\Sp(1)) = \HH{\PP}^1 \approx S^4$. As the
first of the $\Sp(1)$ factors in $H$ acts trivially on the fibre $\C^2=\HH$,
the bundle $G \times_H \HH$ is the negative spinor bundle. The Yang-Mills
connection (with constant norm) that is used in the bundle construction is
again the one compatible with the self-dual positive Einstein structure on
$S^4$. The $\Spin(7)$ holonomy metric is explicit and is asymptotic to the
Ricci-flat cone over the constant curvature metric on $S^7$.

By writing $S^7$ as the homogeneous space $G/K=(\Sp(2)\times \U(1))/\Sp(1)(\Delta \U(1))$
(and thereby increasing the number of parameters of homogeneous metrics from
$2$ to $3$),  Cveti{\u{c}}-Gibbons-L\"u-Pope constructed a $1$-parameter family of
complete $\Spin(7)$-metrics on the negative spinor bundle over $S^4$ \cite{CGLP3}
(see Theorem \ref{nonnullthm}(4) and \S {\bf 5C (4)}).
These metrics are explicit in the sense that the metric components can be
expressed in term of hypergeometric functions, and they represent solutions by
quadratures of the $\Spin(7)$ ODEs. Unlike the metrics in the previous paragraph,
these have ALC asymptotics (the base of the cone is $\C\PP^3$ with the Ziller
metric \cite{Zi}).

Numerical evidence for another $1$-parameter family of complete ALC $\Spin(7)$
metrics was given in \cite{CGLP4}.

\smallskip

$\Spin(7)$ holonomy metrics having the Aloff-Wallach spaces $G/K=\SU(3)/\U_{p,q}$ as
principal orbit type have been intensely studied. Except for one situation below,
we shall assume that $(p, q)=1$. We first point out two simple facts which, if
not kept in mind, would possibly lead to misconceptions. Although the coset spaces
associated to $\{p, q, m=-p-q\}$, $\{-p, -q, p+q\}$, or any permutation of these
triples are equivariantly diffeomorphic, there {\em is} a difference when
one works with a specific singular orbit $G/H$. For example, if we choose $H=S(\U(2)\U(1))
\subset \SU(3)$ (so that the singular orbit is $\C\PP^2$), then, as was pointed out
already in \cite{CGLP4}, $H/K$ has fundamental group $\Z/|p+q|\Z$, so smoothness
requires us to take $p+q = \pm 1$. In particular, the subgroups $\U_{1, -2}$
and $\U_{2, -1}$ are allowed but not $\U_{1,1}$.

A second point is that all the Aloff-Wallach spaces, including $\SU(3)/\U_{1,-1}$,
admit, up to the action of the gauge group, two $\SU(3)$-invariant Einstein metrics
(see the discussion in Example \ref{Wallach1}). However, in the physics literature,
e.g., \cite{CGLP4}, \S 3.1, it is sometimes suggested that $\SU(3)/\U_{1,-1}$ has only
one $\SU(3)$-invariant Einstein metric. The difference between the two exceptional
Aloff-Wallach spaces $\SU(3)/\U_{1,1}$ and $\SU(3)/\U_{1,-1}$ is that in the former
case, all the $\SU(3)$-invariant Einstein metrics have a representative that is
diagonal with respect to the usual decomposition of the isotropy representation
(see Example \ref{Wallach2}). In the latter case, this is false because the gauge
orbit of one of the Einstein metrics lies inside the non-diagonal invariant metrics.
The consequence for cohomogeneity one metrics is that if one uses a $5$-functions
ansatz, as is done in \cite{KaY} and \cite{Re3}, then one cannot see the Ricci-flat
cone on the non-diagonal Einstein metric in the asymptotics of the $\Spin(7)$ metrics
being studied.

$\Spin(7)$ holonomy metrics of cohomogeneity one based on the Aloff-Wallach spaces
were first studied by various teams of physicists, e.g., \cite{CGLP1}, \cite{CGLP3},
\cite{CGLP4}, \cite{GuSp}, \cite{KaY}, because of their relevance in M-theory.
Now physicists are not just interested in the metrics, but also in branes
(aka calibrated submanifolds), string dynamics, and other issues in quantum
gravity. Even though relatively few rigorous mathematical results were obtained,
many important issues and phenomena were brought forth by these works, as well as
numerical evidence for various families of special holonomy metrics and their
relation to each other. Still, there are some analytic results.

In \cite{CGLP4}, a special exact solution to the $\Spin(7)$ equations was obtained
for the generic Aloff-Wallach space $\SU(3)/\U_{p,q}$ whenever $2p > q \geq 0$ holds.
The singular orbit is $\C\PP^2$ and the metric has ALC asymptotics (so is complete at
infinity). Unfortunately, this metric is not smooth--there is a conical singularity
in each normal space to the singular orbit. This solution still makes sense when
$p=1, q=0$ (one in effect looks at the subset of cohomogeneity one metrics which are
diagonal on the principal orbit), and one obtains a smooth $\Spin(7)$ metric in this case,
which was independently discovered in \cite{GuSp}.

The physicists paid special attention to the principal orbits $\SU(3)/\U_{1,1}$
and the diffeomorphic spaces $\SU(3)/\U_{1,-2}$ and  $\SU(3)/\U_{-2,1}$. In the
last two situations, taking $H={\rm S}(\U(2)\U(1))$ as the singular isotropy group,
one gets a smooth vector bundle over the corresponding $\C\PP^2$. In the first case,
one cannot have smooth cohomogeneity one metrics near the singular orbit since
$H/K$ is a real projective space. However, it was realised in \cite{KaY} that one
can replace the principal orbit by a suitable $\Z/2\Z$ quotient for which the new
$H/K$ is a sphere. It turns out that $(\SU(3)/\U_{1,1})/\Z_2$ has a $3$-Sasakian
structure, and, exploiting this, Bazaikin was able to construct a $2$-parameter family
of complete $\Spin(7)$ metrics with ALC asymptotics \cite{Baz1} on a certain complex
line bundle over the coadjoint orbit $\SU(3)/T^2$. Evidence for this family was given
in \cite{KaY} by a perturbative argument.

A second $2$-parameter family of complete $\Spin(7)$ metrics was described in
\cite{Baz2}. These again have principal orbit $(\SU(3)/\U_{1,1})/\Z_2$, but the
singular orbit is now $\C\PP^2$. The metrics have ALC asymptotics in general
but when the parameters coincide, then they have AC asymptotics. Evidence for
this family was again given in \cite{KaY}.

Notice that in \cite{KaY} the same principal orbit type gave rise to an explicit
$1$-parameter family of Calabi-Yau ($\SU(4)$ holonomy) metrics on the same complex
line bundle over $\SU(3)/T^2$ with AC asymptotics. At the limiting value of the
parameter, the topology of the bundle jumps and one obtains Calabi's hyperk\"ahler
metrics on $T^*\C\PP^2$. See also \cite{BazM}.

Recently, Reidegeld has carried out a systematic study of the construction of
$\Spin(7)$ holonomy metrics of cohomogeneity one. The first question addressed by him is:
which principal orbits $G/K$ are admissible ? It is known that they have to have a
$G$-invariant ${\rm G}_2$ structure. The classification of such principal orbits
is given in \cite{Re1}. If the invariant ${\rm G}_2$ structure is in addition
cocalibrated (i.e., its Hodge star is closed), then a theorem of Hitchin \cite{Hi2}
yields, via the solution to a flow equation, a $\Spin(7)$ holonomy metric on $I \times (G/K)$
where $I$ is some open interval. In other words, Hitchin's flow equation specialized
to homogeneous cocalibrated ${\rm G}_2$ structures is equivalent to the first order
ODE system expressing the $\Spin(7)$ condition, which is often deduced in the physics
literature via a superpotential (cf \S 5).

In \cite{Re3}, the case of the Aloff-Wallach spaces as principal orbits is carefully
analysed.  Using both the cocalibrated condition as well as the $\Spin(7)$ holonomy
condition, Reidegeld  determines all the possible singular orbits. These are
$\SU(3)/T^2$ and $\C\PP^2$ in general, and for the cases $\SU(3)/\U_{1, -1},
\SU(3)/\U_{1, 0},$ or $\SU(3)/U_{0, 1}$, we can also have $\SU(3)/\SU(2) \approx S^5$ and
the isotropy irreducible homology sphere $\SU(3)/\SO(3)$. For each of these possibilities,
the local existence (i.e., existence in a tubular neighbourhood of the singular orbit)
of general cohomogeneity one Einstein metrics as well as $\Spin(7)$ holonomy metrics is
examined in detail by using and adapting the methods in \cite{EW}. For example, in
the case of a generic Aloff-Wallach space, Reidegeld shows that only $\C\PP^2$ can
occur as a singular orbit for $\Spin(7)$ holonomy metrics, and in this case there is
a maximal $2$-parameter family of local solutions. (Note that smoothness requires $|p+q|=1$.)
The numerical evidence for this family of metrics was given in \cite{CGLP4}.
Recall that the Painl{\'e}ve analysis \cite{DW7} gives a $4$-parameter family of
solutions (two of the parameters are trivial, corresponding to homothety and the
location of the singularity). These represent $\Spin(7)$ holonomy metrics on the infinite end
and have ALC asymptotics. The global existence question, however, seems to be still
open. We refer the reader to \cite{Re3} for the detailed results for the exceptional
Aloff-Wallach spaces, including a uniqueness theorem \cite{Re4} for the case where
the principal orbit is $\SU(3)/\U_{1,1}$ and the singular orbit is $\SU(3)/T^2$.

\begin{rmk} \label{newreidegeld}
An analysis along the same lines for cohomogeneity one metrics with holonomy lying
in $\Spin(7)$ and for which the group $G$ is $\SU(2)\times \SU(2) \times \SU(2)$
or $ \SU(3) \times \SU(2)$ is given in \cite{Re2}. In particular in each case there
is only one possible principal orbit, and its isotropy representation has no
multiplicities. The principal orbits are circle bundles over $S^2 \times S^2 \times S^2$
and $\C\PP^2 \times \C\PP^1$ respectively. The holonomy of the possible cohomogeneity
one metrics turn out all to be $\SU(4)$, and so the smooth complete metrics which occur
in \cite{Re2} should be among those given in \cite{DW1}.
\end{rmk}

\end{document}